\documentclass[onefignum,onetabnum]{siamart190516}

\usepackage{lipsum}
\usepackage{amsfonts}
\usepackage{amssymb}
\usepackage{amsmath}
\usepackage{graphicx}
\usepackage{tikz}
\usepackage{xcolor}
\usepackage{epstopdf}
\usepackage{algorithmic}
\usepackage{xr}
\usepackage{pgf}
\usepackage{pgfplots}
\usepackage{MnSymbol}
\usepackage{subfig}

\usetikzlibrary{intersections}
\usetikzlibrary{calc}
\usetikzlibrary{shapes}

\ifpdf
  \DeclareGraphicsExtensions{.eps,.pdf,.png,.jpg}
\else
  \DeclareGraphicsExtensions{.eps}
\fi

\newsiamremark{remark}{Remark}
\newsiamremark{hypothesis}{Hypothesis}
\crefname{hypothesis}{Hypothesis}{Hypotheses}
\newsiamthm{claim}{Claim}

\headers{Intentionally left blank}{Intentionally left blank}

\title{An FMM accelerated Poisson Solver for
Complicated Geometries in the Plane Using Function Extension\thanks{\funding{The first author gratefully acknowledges the support from the Knut and Alice Wallenberg Foundation under grant 2020.0258.}}}
\author{Fredrik Fryklund\thanks{Courant Institute of Mathematical Sciences, New York University,
New York, NY 10012, United States
  (\email{fredrik.fryklund@cims.nyu.edu}, \email{greengard@cims.nyu.edu}).}
\and Leslie Greengard\footnotemark[2]\textsuperscript{ ,}\thanks{Center for Computational Mathematics, Flatiron Institute, Simons Foundation, New York, NY 10010, United States
(\email{lgreengard@flatironinstitute.org}).}}
\usepackage{amsopn}

\makeatletter
\newcommand*{\addFileDependency}[1]{
  \typeout{(#1)}
  \@addtofilelist{#1}
  \IfFileExists{#1}{}{\typeout{No file #1.}}
}
\makeatother

\newcommand*{\myexternaldocument}[1]{%
    \externaldocument{#1}%
    \addFileDependency{#1.tex}%
    \addFileDependency{#1.aux}%
}

\newcommand{\bPhi}{{\mathbf \Phi}}
\newcommand{\bA}{{\mathbf A}}
\newcommand{\Chi}{{\mathcal P}}
\newcommand{\bbp}{{\mathbf{p}}}

\ifpdf
\hypersetup{
  pdftitle={draft},
  pdfauthor={F. Fryklund, L. Greengrad}
}
\fi

\myexternaldocument{supplement}
\headers{Adaptive Poisson Solver in Complicated Domains}{F. Fryklund and L. 
Greengard}

\begin{document}

\maketitle

\begin{abstract}
We describe a new, adaptive solver for the two-dimensional 
Poisson equation in complicated 
geometries. Using classical potential theory, we represent the 
solution as the sum of a volume potential and a double layer potential. 
Rather than evaluating the volume potential over the given domain,
we first extend the source data to a geometrically simpler region
with high order accuracy. This allows us to accelerate the evaluation 
of the volume potential using an efficient,
geometry-unaware fast multipole-based algorithm.
To impose the desired boundary condition, it remains only
to solve the Laplace equation with suitably modified boundary data. 
This is accomplished with existing fast and accurate boundary integral 
methods. The novelty of our solver is the scheme used for creating the 
source extension, assuming it is provided on an adaptive quad-tree. 
For leaf boxes intersected by the boundary, 
we construct a universal ``stencil" and
require that the data be provided at the subset of those points 
that lie within the domain interior.
This universality permits us to  
precompute and store an interpolation matrix which
is used to extrapolate the source data
to an extended set of leaf nodes with full tensor-product grids on each.
We demonstrate the method's speed, robustness and high-order convergence 
with several examples, including domains with piecewise smooth boundaries.
\end{abstract}



\section{Introduction}
\label{sec:intro}
We consider the problem of rapidly and accurately solving the Poisson equation 
\begin{align}
 \label{eq:poissoneq}
\Delta u(\mathbf{x}) &= f(\mathbf{x}), \quad \text{for }\mathbf{x} \in \Omega,\\
   \label{eq:poissoneqbc}
u(\mathbf{x}) &= g(\mathbf{x}), \quad \,\text{for } \mathbf{x} \in \partial\Omega,
\end{align}
in complicated domains in the plane. 
Here, $u$ is an unknown function, $f$ is a smooth source density and  
$g$ is the specified Dirichlet boundary data.
While many fast solvers are based on direct discretization of the partial
differential equation itself, recent
years have witnessed substantial progress in developing solvers based
on potential theory, that make use of the linearity of the problem
to solve \eqref{eq:poissoneq}, \eqref{eq:poissoneqbc} in two steps.
One first computes  a ``particular solution" $v(\mathbf{x})$
that satisfies 
\begin{equation}
 \label{eq:poissoneqf}
\Delta v(\mathbf{x}) = f(\mathbf{x}) \quad \text{for }\mathbf{x} \in \Omega,
\end{equation}
without regard to the boundary condition,
and then finds a harmonic function $w(\mathbf{x})$ that satisfies
\begin{align}
 \label{eq:laplaceeq}
\Delta w(\mathbf{x}) &= 0 \quad \text{for }\mathbf{x} \in \Omega,\\
w(\mathbf{x}) &= g(\mathbf{x})- v(\mathbf{x}), \quad \,\text{for } \mathbf{x} \in \partial\Omega. \nonumber
\end{align}
Clearly, $u(\mathbf{x})= w(\mathbf{x}) + v(\mathbf{x})$ is the 
desired solution.
Strong arguments for this approach are that 
\eqref{eq:poissoneqf} can be solved by an integral transform without any
volumetric unknowns and  that
\eqref{eq:laplaceeq} can be solved using a boundary integral equation with
unknowns only on the surface $\partial \Omega$ 
(see, for example, \cite{ASKHAM20171,fryklund2020integral,MGM}).

One possible choice for $v(\mathbf{x})$ is the {\em volume potential}
\begin{equation}
\label{eq:volpot}
\mathcal{V}_{\Omega}[f](\mathbf{x}) = \int\limits_{\Omega}G(\mathbf{x}-\mathbf{y})f(\mathbf{y})\,\mathrm{d}\mathbf{y}, \quad \,\text{for } \mathbf{x} \in \Omega,
\end{equation}
where $G(\mathbf{x})$ is the free-space Green's function 
\cite{evans1998partial}, given in the two-dimensional case by
\begin{equation}
\label{eq:greenf}
  G(\mathbf{x},\mathbf{y}) = \frac{1}{2\pi}\log\|\mathbf{x} - \mathbf{y}\|.
\end{equation}

\begin{definition}
When $\Omega$ is a square with
$f(\mathbf{x})$ given at tensor product grid points on the leaf nodes
of an adaptive quad-tree data structure,
highly optimized fast multipole methods are available
for computing volume potentials of the form \eqref{eq:volpot}
\cite{vfmm,malhotra_biros_2015}. 
We will refer to such methods as {\em volume-integral 
fast multipole methods (VFMMs)}.
\end{definition}

\begin{remark}
VFMMs assume that $f(\mathbf{x})$ is resolved with high order accuracy
by a piecewise polynomial approximation on the leaf nodes.
For orders of accuracy greater than four, VFMMs 
typically use tensor product Chebyshev or Legendre grids on the 
leaf nodes for stable high order approximation. 
For those familiar with VFMMs, recall that
\eqref{eq:volpot} is computed exactly (for the piecewise polynomial
approximation of the source density) in the near field and with 
arbitrary, user-controlled precision in the far field.
\end{remark}

For general domains, however, VFMMs cannot be applied directly with 
high order accuracy,
since there will be {\em cut leaf nodes} that are intersected by the
boundary $\partial \Omega$ and where the data is only defined in the 
domain interior. This prevents simple high order polynomial approximation
of $f(\mathbf{x})$ on the cut nodes, since the function is not
locally smooth. In this paper, we seek to 
enable the application of VFMMs by first extending 
the function $f(\mathbf{x})$ smoothly to a function 
$f^e(\mathbf{x})$ defined on a region $E \supset \Omega$ 
for which a VFMM can be applied (see Fig. \ref{fig:makext}).

The combination of function extension and fast solvers is an active area
of research.
In \cite{STEIN2016252,STEIN2017155}, for example,
an extension of the source density $f(\mathbf{x})$ is obtained
through an immersed boundary formalism.
In \cite{ASKHAM20171}, function extension is carried out
using a boundary integral formulation
(with harmonic extension yielding
a $C^0$ extension, biharmonic extension yielding
a $C^1$ extension, etc.).  Once the extended function has been obtained, 
the fully adaptive solver of \cite{ASKHAM20171} computes the extended volume 
integral using the VFMM algorithm of \cite{vfmm,greengard1996direct}. 
Another approach is Fourier continuation. One such  
scheme is described in \cite{doi:10.1137/20M1373189}, 
where the source density is extended 
normal to the boundary, through projection onto a basis that vanishes in the 
vicinity of the boundary. 
For a good discussion of Fourier-based extension, see \cite{boydext}.
In the active penalty method \cite{activepenalty},
an extension is created by matching boundary data and normal derivatives 
up to order $k$ in terms of a carefully crafted set of basis functions
which gives an extension with global regularity $C^k$.

Function extension is not the only way in which the computation of volume
potentials can be accelerated using VFMMs. 
One alternative is to modify the VFMM algorithm 
to treat the cut leaf
nodes 
via more elaborate approximation and quadrature tools that depend
on the precise intersection of $\partial \Omega$ with the leaf node
(see \cite{shravan_hai_volint} and the references therein).
Another alternative is the recently developed
technique of {\em function intension} \cite{functionintension}, 
where the source density is coverd by regular tensor product leaf nodes
on an adaptive quad-tree in the interior of $\Omega$, 
blended with a conforming mesh in the neighborhood of the boundary.
This permits the use of a VFMM for the interior degrees of 
freedom, but needs to be coupled to an auxiliary fast Poisson solver on 
a tubular neighborhood of the boundary.

In this paper, we describe our new function extension algorithm
in detail, using the Dirichlet problem for the Poisson equation as
our model. The method is fully adaptive in the interior of the domain, 
high-order accurate, robust, and fast.
The main novelty lies in creating 
a universal, level-independent, oversampled interpolation matrix, 
recruiting sufficient data from the neighbors of cut
leaf nodes, and defining the range
of the extension based on the local mesh size of the adaptive discretization.
Our method is
similar to the {\em partition of unity extension scheme} (PUX)
developed in  \cite{fryklund2018partition} for uniform grids.
In the present scheme, however, blending and partitions of unity are avoided. 
Instead, the extension for each cut leaf node is entirely local,
based only on data from the square itself or its nearest
neighbors. Moreover, there is no need to 
truncate the function smoothly to zero; 
we simply extend it to cover a domain $E \supset \Omega$ which is discretized
as a collection of leaf nodes beyond which 
$f(\mathbf{x})$ is identically zero (Fig. \ref{fig:makext}).
(In PUX, the data is represented on a uniform grid, 
rather than a quad-tree, extended smoothly to zero, with a global
interpolation framework based on the fast Fourier transform.)
Rather than \eqref{eq:volpot}, the VFMM then 
computes the volume potential
\begin{equation}
 \label{eq:volinte}
v(\mathbf{x}) =
\mathcal{V}_{E}[f^e](\mathbf{x}) 
= \int_E G(\mathbf{x-y}) f^e(\mathbf{y}) \, d{\mathbf{y}}
\end{equation}
where $f^e$ is the smooth extension of $f$.

An important feature of our method is that it is
agnostic to the smoothness of the boundary. It simply assumes that 
the adaptive quad-tree has resolved the source data well enough, and that
the user can identify points as being either inside or outside the domain.
We will demonstrate that with an eighth order accurate VFMM,
we obtain an eighth order accurate scheme for the full problem,
even on piecewise smooth domains.
We will also show that the overall approach is compatible with other 
extension schemes, including one-dimensional extensions along lines, 
using either the rational function approximation of \cite{baryrat} or the 
diffeomorphism-based method of 
\cite{epsteinjiang2022}. For a review of 

There are, of course, drawbacks to extensions schemes - the major ones
being caustics and ill-conditioning. The former arise when a domain 
boundary curves back on itself, so that two exterior normals intersect 
close to the domain. This can be overcome by ensuring that 
the length scale of leaf nodes in the quad-tree near such points
must be commensurate with the distance to the nearest such intersection.
The difficulty is that this constraint could result in excessive 
refinement, even though the geometry and the 
data may be simple to resolve. Ill-conditioning is an inherit concern with 
function extension, since it is an extrapolation process.
This effect is mitigated by the fact that,
as the quad-tree is refined, the data becomes locally smoother on the
scale of the leaf node and the extension problem becomes simpler as well. 
A detailed analysis of the conditioning of the process 
remains to be carried out, but
experiments indicate that our method performs well without excessive
resolution.
The algorithm requires that data be provided at
auxiliary nodes close to the boundary, but this is to be expected
in a high-order formulation, and the node locations are specified
as soon as the quad-tree is created, so can be considered part of the 
discretization process.

This paper is organized as follows. In \cref{sec:mathprel},
we review the needed elements of classical potential theory for the 
Poisson equation, and in \cref{sec:fext}, we discuss function extension with Gaussians. In \cref{sec:volpot}, we present the data structures used to 
discretize the right-hand side $f(\mathbf{x})$ and create 
its extension. Layer potentials are discussed in more detail in 
\cref{sec:dblpot} and the performance of the algorithm is illustrated
in \cref{sec:numres}, along with a discussion of some implementation details. 
In \cref{sec:concl}, we discuss extensions of the present scheme
and consider avenues for future improvement.

\section{Mathematical preliminaries}
\label{sec:mathprel}

Let $D$ be an open, 
bounded subset of $\mathbb{R}^{2}$, 
which is either simply or multiply connected.
For a 
point $\mathbf{x}$ in $\mathbb{R}^{2}$, we will denote its
Cartesian components by $(x_{1},x_{2})$ and its Euclidean norm by 
$\| \mathbf{x} \|$. 
For $\mathbf{x}, \mathbf{y}$ in $\mathbb{R}^{2}$, their inner product
will be denoted by $\mathbf{x} \cdot \mathbf{y}$.
Unless otherwise stated, we assume that the domain
has a boundary $\partial \Omega$ which is 
at least twice continuously differentiable.
In the case of an {\em interior} problem, $\Omega = D$
and the problem is fully specified. 
In the case of an \emph{exterior} problem,
$\Omega = \mathbb{R}^2\setminus \bar{D}$, 
in which case we must also specify a condition at infinity for the 
problem to be well-posed. That is, we must specify a constant $A$ such that 
\begin{equation}
\label{radcond}
u(\mathbf{x}) \rightarrow A \log(\| \mathbf{x} \|) 
\end{equation} 
as $\| \mathbf{x} \| \rightarrow \infty$
\cite{FollandPDE,Kelloggpot}. Bounded solutions correspond to setting
$A=0$.

For simplicity, let us begin with
the interior problem in a simply-connected domain.
By standard potential theory \cite{FollandPDE,guenther1988partial,Kelloggpot}, 
an explicit representation of the solution $u$ can be formulated as
\begin{equation}
\label{eq:solution} 
  u(\mathbf{x}) =  \mathcal{V}_{E}[f^e](\mathbf{x}) + \mathcal{D}[\sigma](\mathbf{x}), \quad \,\text{for } \mathbf{x} \in \Omega,
\end{equation}
where the volume potential $\mathcal{V}_{E}[f^e](\mathbf{x})$ is defined in  
\eqref{eq:volinte}, so long as 
$E \supset \Omega$ and $f^e = f$ within $\Omega$. Here,
\begin{equation}
\label{eq:dlpot}
\mathcal{D}[\sigma](\mathbf{x}) = \int\limits_{\partial\Omega}\frac{\partial G(\mathbf{x}-\mathbf{y})}{\partial \boldsymbol{\nu}(\mathbf{y})}\sigma(\mathbf{y})\,\mathrm{d}s, \quad \,\text{for } \mathbf{x} \in \Omega,
\end{equation}
is the \emph{double layer potential}, with unknown 
\emph{layer density} $\sigma\colon \partial \Omega \rightarrow \mathbb{R}$,
$\boldsymbol{\nu}(\mathbf{y})$ denotes
the unit normal at the point $\mathbf{y} \in \partial \Omega$, and
$\frac{\partial G(\mathbf{x}-\mathbf{y})}
{\partial \boldsymbol{\nu}(\mathbf{y})}$ denotes the normal derivative
of the Green's function \eqref{eq:greenf}.
It is straightforward to see that $\mathcal{D}[\sigma](\mathbf{x})$ is
harmonic and that
the kernel of $\mathcal{D}[\sigma]$ in \cref{eq:dlpot} is 
\begin{equation}
\label{eq:dblker}
  \frac{\partial G(\mathbf{x}-\mathbf{y})}{\partial \boldsymbol{\nu}(\mathbf{y})} = \frac{1}{2\pi}\frac{(\mathbf{x}-\mathbf{y})\cdot\boldsymbol{\nu}(\mathbf{y})}{\|\mathbf{x}-\mathbf{y}\|^{2}}.
\end{equation}
Note that the limiting value of \cref{eq:dblker},
as $\mathbf{x}$ approaches $\mathbf{y}$ along the boundary, 
is $-\tfrac{1}{2\pi}\kappa(\mathbf{y})$,
where $\kappa(\mathbf{y})$ is the curvature at $\mathbf{y}$. Thus, 
assuming the boundary is at least twice differentiable, the kernel is a 
continuous function. For a $C^{k+2}$ boundary, the kernel is $C^k$.

In order to satisfy the desired Dirichlet boundary conditions,
\cref{eq:poissoneqbc}, we seek a layer density $\sigma$ such that 
$\mathcal{D}[\sigma](\mathbf{x}) = g(\mathbf{x}) - \mathcal{V}_{E}[f^{e}](\mathbf{x})$ for $\mathbf{x}$ on $\partial\Omega$. 
This is achieved by taking the limit of \eqref{eq:solution} 
as $\mathbf{x}$ approaches the boundary from the interior and applying
standard jump conditions \cite{KressRainer2014LIE}, yielding the integral 
equation
\begin{equation}
\label{eq:inteq}
  (-\tfrac{1}{2} + \mathcal{D})[\sigma](\mathbf{x}) = g(\mathbf{x}) - \mathcal{V}_{E}[f^{e}](\mathbf{x}),\quad \text{for } \mathbf{x} \in \partial\Omega.
\end{equation}
\Cref{eq:inteq} is a Fredholm integral equation of the second kind 
for $\sigma$, since $\mathcal{D}[\sigma]$ is a compact operator 
with a continuous kernel on a $C^2$  boundary (as noted above).
It follows by the Fredholm alternative that \cref{eq:inteq} has a unique 
solution \cite{atkinson_1997}.
Once $\sigma$ has been obtained, we have a complete solution to the 
full problem.

\begin{remark}
For Neumann boundary value problems, where \eqref{eq:poissoneqbc}
is replaced by
\begin{equation}
   \label{eq:poissoneqnbc}
\frac{\partial u}{\partial \boldsymbol{\nu}} (\mathbf{x}) = 
g(\mathbf{x}), \quad \,\text{for } \mathbf{x} \in \partial\Omega,
\end{equation}
the approach is essentially the same, except that the homogeneous solution 
is expressed as a \emph{single layer potential}
\begin{equation}
\label{eq:slpot}
\mathcal{S}[\sigma](\mathbf{x}) = \int\limits_{\partial\Omega} 
G(\mathbf{x}-\mathbf{y}) \, \sigma(\mathbf{y})\,\mathrm{d}s, \quad 
\,\text{for } \mathbf{x} \in \Omega.
\end{equation}
Imposing  \eqref{eq:poissoneqnbc}
leads to a second kind Fredholm equation, ensuring a unique solution 
(up to an arbitrary constant) so long as 
$\int_{\partial\Omega} g(\mathbf{x}) \, ds = 0$.
\end{remark}

For the exterior Dirichlet problem in 
$\Omega = \mathbb{R}^2\setminus \bar{D}$, 
we first compute a particular solution
of the form \eqref{eq:volinte}, where the extension is now {\em into} $D$.
The exterior harmonic correction is then represented in the form
\begin{equation}
\label{eq:dlpotext}
u^{H}(\mathbf{x}) =
 \mathcal{D})[\sigma](\mathbf{x}) + \frac{1}{2 \pi} \,
\int\limits_{\partial\Omega}
\sigma(\mathbf{y})\,\mathrm{d}s \, + \,  
\alpha \, \log\|\mathbf{x} - \mathbf{x}_{D}\|,
\end{equation}
where $\mathbf{x}_{D}$ lies in $D$.
Letting $Q = \frac{1}{2 \pi} 
\iint_{E}f^{e}(\mathbf{y})\,\mathrm{d}\mathbf{y}$, we impose
the additional constraint
\begin{equation}  
\alpha = A - Q
\end{equation}  
to ensure the user-specified radiation condition  
\eqref{radcond}. For a discussion of uniqueness of the resulting integral
equation, see \cite{GREENBAUM1993267,mikhlin}. 

\subsection{Multiply-connected domains}

We now consider the interior problem for a multiply connected domain, 
whose boundary consists of $(N_{\Omega} + 1)$ closed curves. 
The outer boundary curve is denoted $\partial\Omega_{0}$, and the interior 
boundary curves are denoted by 
$\partial\Omega_{1},\ldots,\partial\Omega_{N_{\Omega}}$ 
(see \cref{fig:multdomain}). In this setting, it turns out that
there are $N_{\Omega}$ nontrivial homogeneous solutions to the boundary 
integral equation \cref{eq:inteq} \cite{FollandPDE}. 
In order to ensure uniqueness, we proceed as 
in \cite{GREENBAUM1993267}, and write the full solution to the
Poisson equation in the form
\begin{equation}
\label{eq:solutionmult}  
  u(\mathbf{x}) =  \mathcal{V}_{E}[f^e](\mathbf{x}) + 
\mathcal{D}[\sigma](\mathbf{x}) + \sum\limits_{k = 1}^{N_{\Omega}}A_{k}\log\|\mathbf{x} - \mathbf{s}_{k}\|, \quad \,\text{for } \mathbf{x} \in \Omega,
\end{equation}
where $\mathbf{s}_{k}$ is a point inside the interior curve $\partial\Omega_{k}$ and $\{A_{k}\}_{k=1}^{N_{\Omega}}$ are unknown constants, with the additional constraints
\begin{equation}
  \label{eq:densconst}
  \int\limits_{\partial\Omega_{k}} \sigma(\mathbf{y})\,\mathrm{d}\mathbf{y} = 0,\quad\, k = 1,\ldots,N_{\Omega}.
\end{equation}
Imposing the Dirichlet boundary conditions together with \cref{eq:densconst} 
leads to an invertible Fredholm equation of the second kind
for the unknowns $\sigma$ and $\{A_{k}\}_{k=1}^{N_{\Omega}}$.

Finally, we consider the Dirichlet problem posed in 
the region exterior to a collection of $N_{\Omega}$ closed curves 
$\partial\Omega_{1},\ldots,\partial\Omega_{N_{\Omega}}$.
It is shown in  \cite{GREENBAUM1993267}, that the representation
\begin{equation}
\label{eq:solutionmulte}  
  u^H(\mathbf{x}) =  
\mathcal{D}[\sigma](\mathbf{x}) + \frac{1}{2 \pi}
\int\limits_{\partial\Omega} \sigma(\mathbf{x}) \, 
\mathrm{d}\mathbf{x} \ +
\sum\limits_{k = 1}^{N_{\Omega}}A_{k}\log\|\mathbf{x} - \mathbf{s}_{k}\|, \quad \,\text{for } \mathbf{x} \in \Omega,
\end{equation}
together with the constraints
\begin{equation}
  \label{eq:densconstext}
\int\limits_{\partial\Omega_{k}} \sigma(\mathbf{y})\,\mathrm{d}\mathbf{y} = 0,\quad\, k = 1,\ldots,N_{\Omega}-1 \quad , \quad
\sum_{k=1}^{N_{\Omega}} A_{k}= A,
\end{equation}
leads to a well-conditioned 
Fredholm equation of the second kind
for $\sigma$ and $\{A_{k}\}_{k=1}^{N_{\Omega}}$.

\begin{remark}
We will also consider domains with piecewise smooth boundaries. 
For such geometries, the double layer operator is no longer 
compact, but there is an extensive literature on the invertibility
of the corresponding integral equation (see \cite{verchota,cjs1989})
and the design of high order methods for its solution 
(see, for example, \cite{brs,rciphelsing,hrs2019}).
\end{remark}

\newcommand{\normtangoutwards}[3]{
    \draw[->, black] (#2) -- ($(#1)!(#2)!(#3)!-4!(#2)$);
}
\newcommand{\normtang}[3]{
    \draw[->, black] (#2) -- ($(#1)!(#2)!(#3)!.5!(#2)$);
}
\begin{figure}
\centering
\begin{tikzpicture}[thick,scale=0.4, every node/.style={scale=0.6}]

\path[]
    (0.1,4.7)   coordinate (A0)
    (2.8,6.9)   coordinate (B0)
    (4,6.8)   coordinate (C0)
    (6,7.5)   coordinate (D0)
    (8,7)   coordinate (E0)
    (8.5,4.5)   coordinate (F0)
    (10,2)   coordinate (G0)
    (7,0.5)   coordinate (H0)
    (5,1)   coordinate (I0)
    (4,0.1)   coordinate (J0)
;

\draw[inner color = gray!60, outer color = gray!60] plot [smooth cycle] coordinates {(A0) (B0) (C0) (D0) (E0) (F0) (G0) (H0) (I0) (J0)};
\normtangoutwards{B0}{C0}{D0}

\path[font={\tiny}]
  (5,6)   coordinate (A1)
  (6.4,7)   coordinate (B1)
  (7.5,5)   coordinate (C1)
  (5.5,5)   coordinate (D1)
;
\draw[inner color = white, outer color = white] plot [smooth cycle, tension=0.8] coordinates {(A1) (B1) (C1) (D1)};

\normtang{B1}{C1}{D1}

\path[font={\tiny}]
  (1.4,4.8)   coordinate (A2)
  (4,6)   coordinate (B2)
  (3.3,4)   coordinate (C2)
;
\draw[inner color = white, outer color = white] plot [smooth cycle, tension=0.8] coordinates {(A2) (B2) (C2)};

\path[font={\tiny}]
  (5.2,2)   coordinate (A3)
  (7.8,3)   coordinate (B3)
  (8.2,1.2)   coordinate (C3)
;
\draw[inner color = white, outer color = white] plot [smooth cycle, tension=0.8] coordinates {(A3) (B3) (C3)};

\draw[dotted,scale=0.3]  ([xshift=60pt,yshift=60pt]C2) -- ([xshift=20pt,yshift=60pt]A3);
     \node[scale=1.6] at (6.2,3.8) {$\Omega$};
     \node[scale=1.6] at (9.7,4.8) {$\partial\Omega_{0}$};
     \node[scale=1.6] at (6,6) {$\partial\Omega_{1}$};
     \node[scale=1.6] at (2.7,4.8) {$\partial\Omega_{2}$};
     \node[scale=1.6] at (7.2,1.9) {$\partial\Omega_{N_{\Omega}}$};

\end{tikzpicture}
\hspace{2cm}
\begin{tikzpicture}[thick,scale=0.4, every node/.style={scale=0.6}]

\path[]
    (0,4.5)   coordinate (A0)
    (2.8,6.9)   coordinate (B0)
    (4,6.8)   coordinate (C0)
    (6,7.5)   coordinate (D0)
    (8,7)   coordinate (E0)
    (8.5,4.5)   coordinate (F0)
    (10,2)   coordinate (G0)
    (7,0.5)   coordinate (H0)
    (5,1)   coordinate (I0)
    (4,3)   coordinate (J0)
;

\draw[white, inner color = black!60, outer color = black!1] plot  coordinates {(0,0) (0,10) (10,10) (10,0)};

\draw[inner color = white, outer color = white] plot [smooth cycle, tension=0.8] coordinates {(A1) (B1) (C1) (D1)};

\path[font={\tiny}]
  (5,6)   coordinate (A1)
  (6.4,7)   coordinate (B1)
  (7.5,5)   coordinate (C1)
  (5.5,5)   coordinate (D1)
;
\draw[inner color = white, outer color = white] plot [smooth cycle, tension=0.8] coordinates {(A1) (B1) (C1) (D1)};

\normtang{B1}{C1}{D1}

\path[font={\tiny}]
  (1.4,4.8)   coordinate (A2)
  (4,6)   coordinate (B2)
  (3.3,4)   coordinate (C2)
;
\draw[inner color = white, outer color = white] plot [smooth cycle, tension=0.8] coordinates {(A2) (B2) (C2)};

\path[font={\tiny}]
  (5.2,2)   coordinate (A3)
  (7.8,3)   coordinate (B3)
  (8.2,1.2)   coordinate (C3)
;
\draw[inner color = white, outer color = white] plot [smooth cycle, tension=0.8] coordinates {(A3) (B3) (C3)};

\draw[dotted,scale=0.3]  ([xshift=60pt,yshift=60pt]C2) -- ([xshift=20pt,yshift=60pt]A3);

     \node[scale=1.6] at (6,6) {$\partial\Omega_{1}$};
     \node[scale=1.6] at (2.7,4.8) {$\partial\Omega_{2}$};
     \node[scale=1.6] at (7.2,1.9) {$\partial\Omega_{N_{\Omega}}$};

\end{tikzpicture}
\label{fig:multdomain}
\caption{Left: An example geometry for the interior problem on an $(N_{\Omega}+1)$ply connected domain. Right: An example geometry for the exterior problem on an $(N_{\Omega})$ply connected domain}
\end{figure}
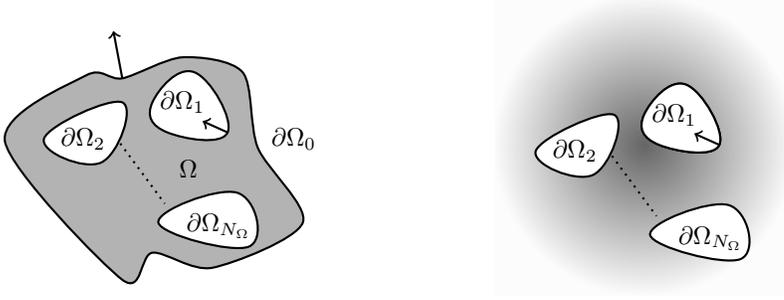
\section{Function extension}
\label{sec:fext}

We turn now to the problem of extending 
the function $f(\mathbf{x})$ defined on $\Omega$ to a function 
$f^e(\mathbf{x}) \in C^{q}(E)$ 
on a region $E \supset \Omega$ 
for which a VFMM can be applied. 
Our scheme is based on local extrapolation using a basis of Gaussians, 
with a precomputed interpolation matrix that can be
obtained using the RBF-QR algorithm\cite{Fornberg2011StableCW},
discussed briefly below. 
This approach is similar to local extension in the PUX algorithm 
\cite{fryklund2018partition}. However, the scheme presented here 
has fewer parameters and requires neither a smooth taper to zero nor 
a blending of multiple local extensions through a partition of unity.

\subsection{Interpolation in a Gaussian basis} \label{sec:gaussinterp}
Consider the approximation
\begin{equation}
  \label{eq:basicinterpolant}
  I_{f}(\mathbf{x}) = \sum\limits_{i = 1}^{N_{\mathcal{G}}}\lambda_{i}\phi_{i}(\mathbf{x}),\quad \mathbf{x} \in \bar{S},
\end{equation}
of a function $f\in C^{q}(\bar{S})$, with $q \geq 0$ on the bounded domain $\bar{S}\subset\mathbb{R}^{d}$ for $d = 1,2,3$, with weights 
$\{\lambda_{i}\}_{i = 1}^{N_{\mathcal{G}}}$. The basis consists of Gaussians
$\phi_{i}(\mathbf{x}) = e^{-\varepsilon^{2}\|\mathbf{x}-\mathbf{c}_{i}\|^{2}}$ 
centered at a set of distinct points 
${\mathcal{G}} = \{\mathbf{c}_{i}\}_{i=1}^{N_{\mathcal{G}}}$ in $\bar{S}$. 
We will refer to $\varepsilon$ as
a \emph{shape parameter}, with smaller values corresponding to flatter
basis functions. 
Clearly, $I_{f} \in C^{\infty}(\bar{S}) \subset C^{q}(\bar{S})$.

Let ${\Chi} = \{{\bbp}_{i}\}^{N_{\Chi}}_{i = 1}$ 
be a set of $N_\Chi$ distinct points in $\bar{S}$ and suppose that we wish
to approximate the function values at ${{\Chi}}$
using the representation \eqref{eq:basicinterpolant}.
The weights $\Lambda = (\lambda_{j})\in\mathbb{R}^{N_{\mathcal{G}}}$ can be 
obtained by solving the linear system
\begin{equation}
  \label{eq:interpeq}
  \bPhi_{\Chi,\mathcal{G}} \Lambda = {\bf f}_\Chi,
\end{equation}
where 
$\bPhi_{\Chi,\mathcal{G}}  \in\mathbb{R}^{N_{\Chi}\times N_{\mathcal{G}}}$ 
with
$\bPhi_{\Chi,\mathcal{G}} [i,j] = (\phi_{j}({\bbp}_{i}))$ 
and ${\bf f}_\Chi = (f({\bbp}_{1}),\dots, f({\bbp}_{N_{\mathcal{G}}}))$. 
If $N_{\Chi}>N_{\mathcal{G}}$, then we solve for $\Lambda$ in a 
least squares sense.

Approximation via a sum of Gaussians is a particular case of 
radial basis function approximation
\cite{LARSSON2005103,Fornberg2011StableCW,schabackrbftheory}, 
and we do not seek to review
the literature here, except to note that 
high order accuracy can be achieved by a careful interplay of 
the shape parameter $\varepsilon$ and $N_{\mathcal{G}}$. 
This requires carefully letting 
$\varepsilon \rightarrow 0$ while 
increasing $N_{\mathcal{G}}$
\cite{LARSSON2005103,Fornberg2011StableCW,schabackrbftheory}. 
If $\varepsilon$ were fixed, convergence would stagnate with $N_{\mathcal{G}}$.
On the other hand, for a fixed $N_{\mathcal{G}}$, the linear system 
\eqref{eq:interpeq} becomes increasingly ill-conditioned as 
$\varepsilon \rightarrow 0$, resulting in oscillatory weights $\Lambda$.
Following \cite{lsrbfqr},
it turns out that one can construct a well-conditioned interpolation problem 
for $\varepsilon \approx 10^{-5}$ on the unit box, 
achieving high order convergence.
This involves reformulating \cref{eq:interpeq} to avoid
explicit use of the weights $\Lambda$.
For this, let 
\[ {\bf f}_{\mathcal{G}}= (f(\mathbf{c}_{1}),\dots,
f(\mathbf{c}_{N_{\mathcal{G}}}))
\in\mathbb{R}^{N_{\mathcal{G}}}. \]
If we formally collocate \cref{eq:basicinterpolant} at $\mathcal{G}$,
then $\Lambda = \Phi_{\mathcal{G},\mathcal{G}}^{-1}{\bf f}_{\mathcal{G}}$ 
and we may rewrite
\cref{eq:interpeq} in the form
\begin{equation}
  \label{eq:defA}
  \bA_{\Chi,\mathcal{G}} {\bf f}_{\mathcal{G}} = {\bf f}_\Chi,
\end{equation}
where $\bA_{\Chi,\mathcal{G}} = 
\bPhi_{\Chi,\mathcal{G}} \Phi_{\mathcal{G},\mathcal{G}}^{-1}$
to directly obtain the desired values ${\bf f}_{\mathcal{G}}$. 
While this formulation avoids $\Lambda$,
it remains to address the ill-conditioning of 
$\Phi_{\mathcal{G},\mathcal{G}}$.
It turns out that stable, accurate solutions can be obtained using the 
RBF-QR method \cite{Fornberg2011StableCW}. The essential 
idea is to expand each Gaussian in an intermediate (well-conditioned) basis 
consisting of a combination of powers, Chebyshev polynomials, and 
trigonometric functions. Leaving out the details, the total cost of
RBF-QR is of the order 
$\mathcal{O}(N_{\Chi}N_{\mathcal{G}}M^{2})$, 
where $M > N_{\mathcal{G}}$ is the number of functions used in the 
intermediate basis. 
This cost would be prohibitive if carried out at every cut
leaf node in our adaptive discretization. However, if the sets 
$\Chi$ and $\mathcal{G}$ are universal, then 
$\bA$ can be precomputed and stored.
In that case, the cost of solving the least squares problem
\eqref{eq:defA} is of the order 
$\mathcal{O}(N_{\Chi}N_{\mathcal{G}}^2 + N_{\mathcal{G}}^3)$.
In the next section,
we describe how to construct such a univeral 
matrix.

\subsection{Extension from cut leaf nodes} \label{sec:leafextend}

\begin{figure}
\centering
\includegraphics[width=3.5in]{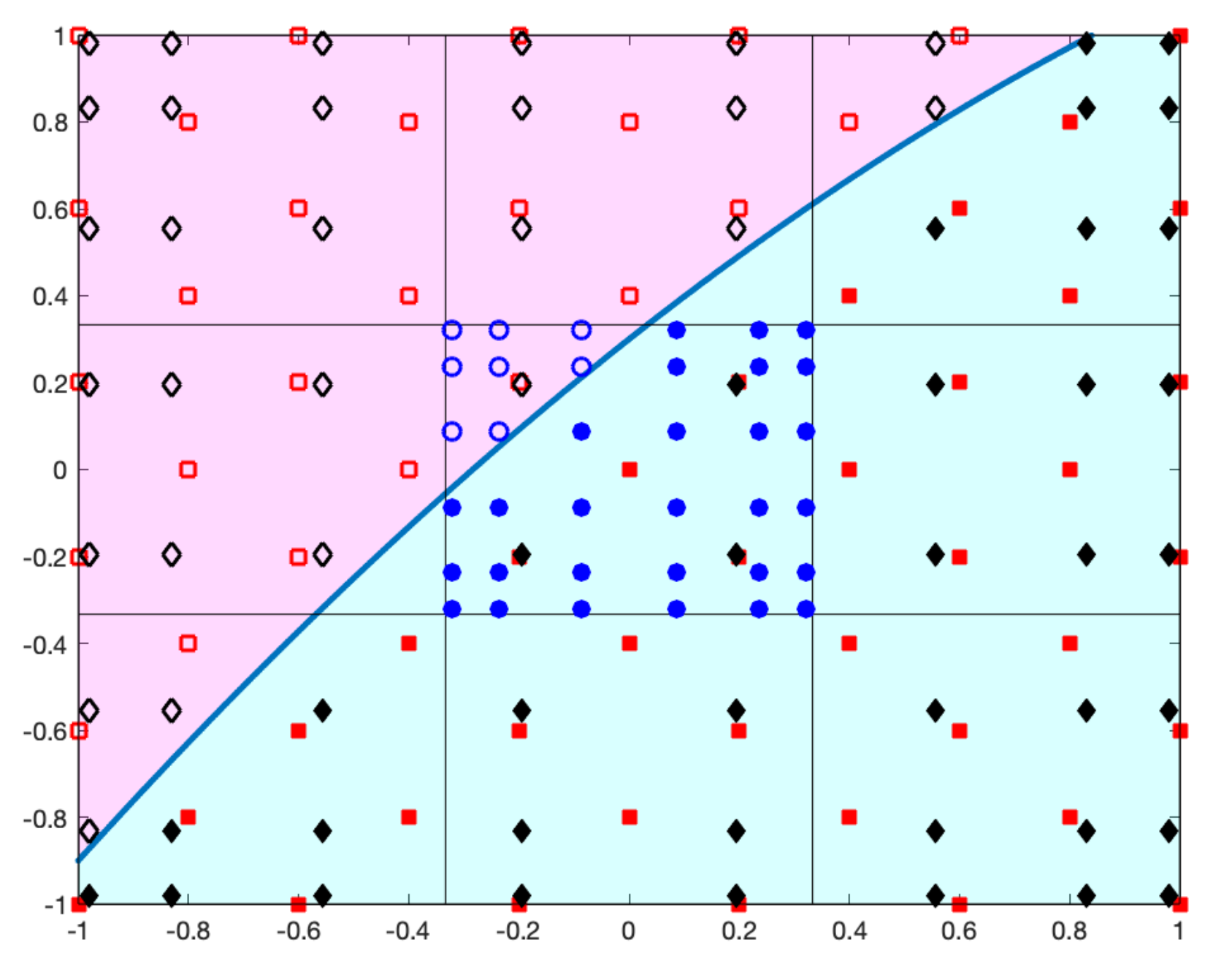}
\caption{A cut cell $S$ of interest (the central square) and its 
8 near neighbors at the same scale. Their union defines the {\em extension 
square} $\bar{S}$.
In this case, suppose that the region below the curve (light blue) 
corresponds to the domain interior and the region above the curve
(light purple) corresponds to the exterior. 
The marked points with a blue, circular shape in $S$ are the 
scaled Chebyshev nodes on $S$. 
The marked points with a black, diamond shape are the scaled 
Chebyshev nodes on $\bar{S}$. 
The marked points with a red, square shape are the {\em support} nodes
for the Gaussian basis functions. 
The points in each of these sets that lie in the domain interior are
indicated with filled markers.
The points in these sets that lie in the domain exterior are
indicated with unfilled markers.
In our extension algorithm, we construct a least-squares interpolant at the 
interior points and use it to obtain values at the 
Chebyshev nodes on $\bar{S}$. 
} 
\label{fig:gridpts}
\end{figure}

Let $S$ be a square of sidelength $L$ which is cut by the boundary $\Gamma$
of our domain $\Omega$,
and let $\mathcal{X}_S$ be the $K \times K$ tensor product Chebyshev grid 
scaled to $S$.
We define the \emph{extension} square $\bar{S}$ to be a square of 
sidelength $3L$, with the same center as $S$ (see Fig. \ref{fig:gridpts}). 
The $\bar{K} \times \bar{K}$ tensor product Chebyshev grid scaled to $\bar{S}$ is 
denoted by $\mathcal{X}_{\bar{S}}$.

On that square, we also impose a
uniform triangulation, and constructing a Delaunay triangulation.
The vertices of that triangulation 
are chosen as the Gaussian support nodes $\mathcal{G}$.
The details of the construction are not so important - 
just that the number be slightly greater than $O(K^2)$ and that they be 
approximately uniformly distributed in the square.
Let $\Chi = \mathcal{X}_S \cup \mathcal{X}_{\bar{S}} \cup
\mathcal{G}$. For any of these point sets, we let the superscript ${\cal I}$ 
refer to the subset that lies in the interior of $\Omega$ and we let
the superscript ${\cal E}$ refer to the subset that lies 
in the exterior of $\Omega$.
Thus,
$\Chi^{\cal I}$ denotes the subset of 
${\Chi}$ that lies in the interior of $\Omega$, and 
$\mathcal{X}_{\bar{S}}^{\cal E}$ denotes the subset of 
${\mathcal{X}_{\bar{S}}}$ that lies in the exterior of $\Omega$.
The full matrix 
$\bA_{\Chi,\mathcal{G}}$ is universal and can clearly be computed and stored.
Extracting the rows corresponding to interior points results in the 
matrix
$\bA_{{\Chi}^{\cal I},\mathcal{G}}$, 
while extracting the rows corresponding
to ${\mathcal{X}_{\bar{S}}^{\cal E}}$ results in
$\bA_{\mathcal{X}_{\bar{S}}^{\cal E},\mathcal{G}}$.
Assuming that the function $f$ is known at 
${\Chi}^{i}$, 
we can obtain its extension $f_{\mathcal{X}_{\bar{S}}^{\cal E}}$ as
\[ 
f_{\mathcal{X}_{\bar{S}}^{\cal E}} =
\bA_{\mathcal{X}_{\bar{S}}^{\cal E},\mathcal{G}}
\bA_{{\Chi}^{\cal I},\mathcal{G}}^{\dagger} 
f_{{\Chi}^{\cal I}}.
\]
where 
$\bA_{{\Chi}^{\cal I},\mathcal{G}}^{\dagger}$ denotes the 
pseudo-inverse of $\bA_{{\Chi}^{\cal I},\mathcal{G}}$.
This yields the missing values
to extend $f$ to a full tensor product Chebyshev grid 
$\mathcal{X}_{\bar{S}}$ on the extension
square ${\bar{S}}$. From this,
we can easily compute $f^e$ at any point in 
$\bar{S}$ by interpolation.

\section{Discretization, data structures, and the volume potential}
\label{sec:volpot}

We turn now to the task of function extension from a
complicated domain $\Omega$ to a larger domain $E$ for which the 
VFMM can be applied with high order accuracy. 
We assume, without loss of generality, that 
$\Omega$ is contained in the unit square $D$ centered at the origin, 
and that the support $E$ of $f^{e}$ (which is only slightly larger than
$\Omega$) is contained in $D$ as well.
We assume that the boundary $\partial\Omega$ is provided in the form
\begin{displaymath}
\label{bdrydef}
  \partial \Omega = \bigcup\limits_{i = 1}^{N_{\Gamma}}\Gamma_{i},
\end{displaymath}
where we refer to the disjoint segments $\{\Gamma_{i}\}_{i=1}^{N_{\Gamma}}$ as \emph{panels}, and each panel is defined by a parametrization
\begin{displaymath}
\label{panelparam}
  \Gamma_{i} = \{\boldsymbol{\gamma}_{i}(t) \in \mathbb{R}^{2} \,\vert\, t \in [-1,1]\}.
\end{displaymath}
We will refer to the length of each panel as 
$|\Gamma_i| = \int_{-1}^1 \| \boldsymbol{\gamma}_{i}'(t) \| \, dt$.

\begin{figure}
\centering
\includegraphics[width=2.5in]{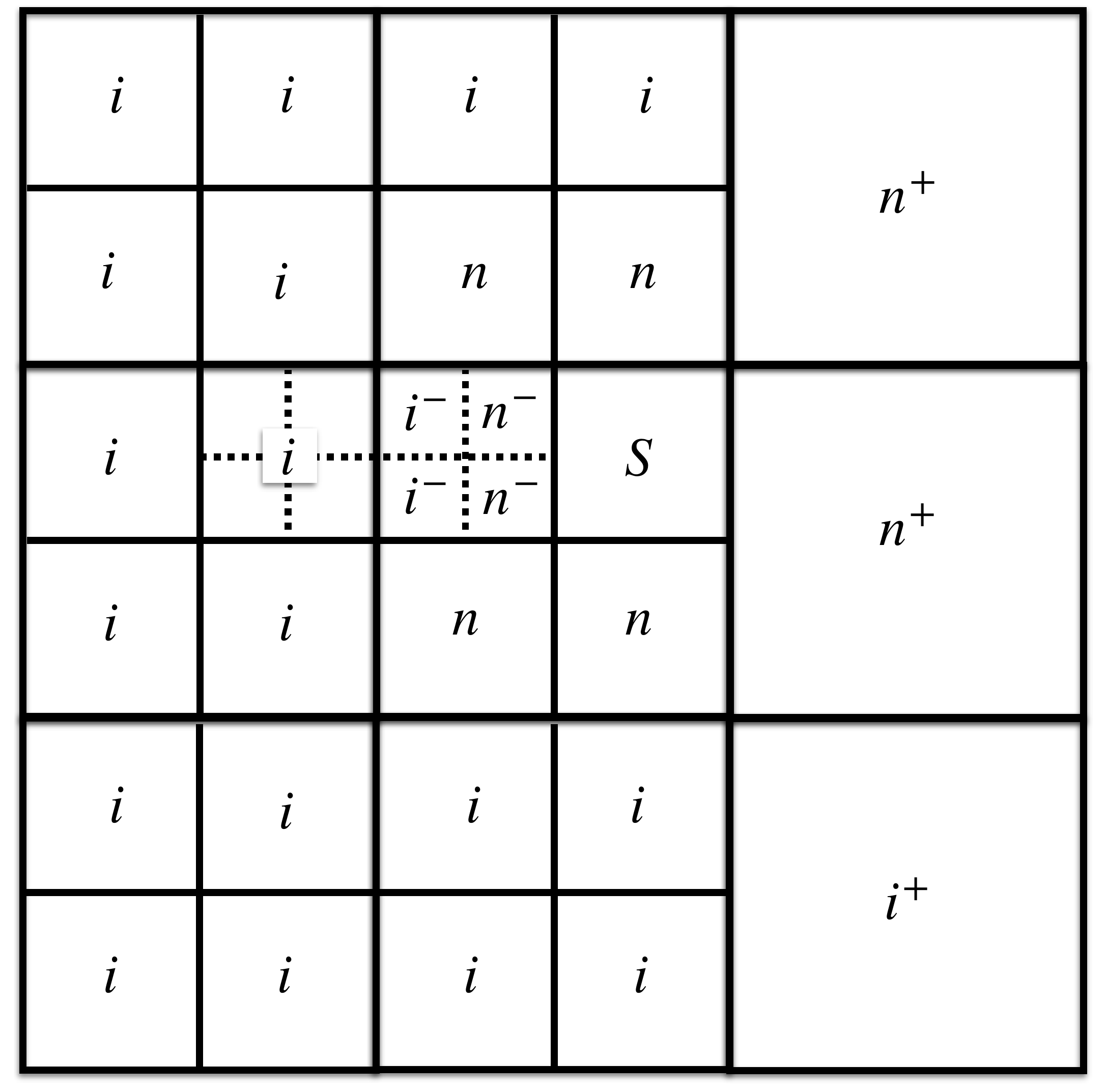}
\caption{An adaptive, level-restricted tree: 
a typical leaf node at level $l$ in the hierarchy (away from the
boundary) is marked
by an $S$. Its colleagues are marked as $n$, its coarse neighbors
as $n^+$ and its fine neighbors as $n^-$. The boxes marked
by $i$ are the children of the colleagues of $S$'s parent
(the so-called {\em interaction list}). 
The box marked $i^+$ is a colleague of $S$'s parent which does not 
touch $S$ and the boxes marked by $i^-$ are children of $S$'s colleagues 
which do not touch $S$. In the VFMM, the 
contributions to $S$ from boxes marked by $i,i^+$, or $i^-$ are
accounted for using multipole and local expansions, while the
contributions from boxes marked by $n,n^+$, or $n^-$ are handled
using precomputed tables. The VFMM is a multi-level algorithm
that computes all such interactions 
using $O(N)$ operations, where $N$ is the total number of points
in the discretization.
} 
\label{fig:tree_levels}
\end{figure}

For the domain itself, 
we assume that an adaptive \emph{quad-tree} is superimposed on 
$D$ to resolve the source density $f(\mathbf{x})$.
For this, the entire box $D$ is referred to as the root node, and
a collection of squares (boxes) at level $l+1$ is obtained by the 
subdivision of some squares (boxes) at level $l$ into four equal parts. 
For a square $S$ at level $l$, the four squares that result from its 
subdivision are referred
to as $S$'s \emph{children}, and $S$ is referred to as their \emph{parent}. 
Squares that do not have children are referred to as \emph{leaf boxes}
or {\em leaf nodes}. For resolve a source distributions with localized
structure, the subdivision process may lead to very fine refinement levels
in some parts of the domain. The only assumption we make about the data
structure is that the tree is {\em level-restricted} or {\em balanced},
meaning that any two leaf nodes which share a boundary point are no more
than one level apart (see \ref{fig:tree_levels}).

\begin{definition}
For a square $S$ at level $l>0$,
its \emph{colleagues} are the boxes at the same refinement level that
share a boundary point with $S$, including itself. 
\emph{Coarse neighbors} of $S$ are leaf nodes at level $l-1$ which 
share a boundary point with $S$ and 
\emph{fine neighbors} are leaf nodes at level $l+1$ which 
share a boundary point with $S$. We define the neighbors of $S$ as the 
union of its colleagues, coarse neighbors and fine neighbors
(\cref{fig:tree_levels}). 
Leaf nodes that lie entirely in the interior of $\Omega$ 
are called \emph{regular leaf nodes}. 
Leaf nodes that are intersected by the boundary are called \emph{cut squares}. 
\end{definition}

For each regular leaf node, we assume that
$f(\mathbf{x})$ is provided on a scaled $K \times K$ tensor product Chebyshev
grid. (In the present paper, we always use Chebyshev nodes of the first kind,
which exclude the endpoints, and fix $K=8$.)
For each cut square $S$ with side length $L$, we define 
the \emph{extension} square $\bar{S}$ as above: the square of length $3L$, 
centered on $S$ (see \cref{fig:gridpts}). 
The extension square can be decomposed into two disjoint subsets: 
the \emph{interpolation} region $\bar{S}^{I}$ that is the intersection of 
$\bar{S}$ and $\Omega$, and the \emph{extension} region 
$\bar{S}^{E} = \bar{S} \setminus \bar{S}^{I}$. 
We define the \emph{extension list} for a cut square $S$ to be
the set of all leaf squares $S'$ intersected by the extension region 
$\bar{S}^{E}$, for which the center of $S$ is the closest of all 
cut squares centers. If two cut square centers are equidistant, 
the latter cut square which has added $S'$ to its extension
list takes precedence.
On each cut square, there is a $K \times K$ Chebyshev grid, for which some 
nodes are within the domain and some not.
On each extension square, we assume there is 
also a $\bar{K}\times \bar{K}$ Chebyshev grid and a set of $N_{\mathcal{G}}$ 
distinct points $C_{\mathcal{G}}$. 

In adaptive refinement, a standard criterion for
regular (non-cut) squares is that the source distribution is resolved.
From tensor product Chebyshev samples, one measure of resolution
is spectral decay: that is, one computes the 
Chebyshev expansion of $f(\mathbf{x})$ and 
requires that the relative
$\ell^{p}$ norm of the vector of Chebyshev coefficients of total order $N-1$ 
be below a prescribed tolerance.
If that is satisfies, the refinement is terminated. Otherwise, one preceeds
to the next level. 

\begin{remark}
In practice, it is simplest to refine a quad-tree without regard to 
level-restriction, based on resolution considerations alone. There are
standard algorithms that take a general adaptive quad-tree as input 
and create a slightly more refined tree which does satisfy the 
level-restriction (see, for example, \cite{biros2008sisc}).
\end{remark}

\begin{figure}[!h]
  \centering
\includegraphics[width=2in]{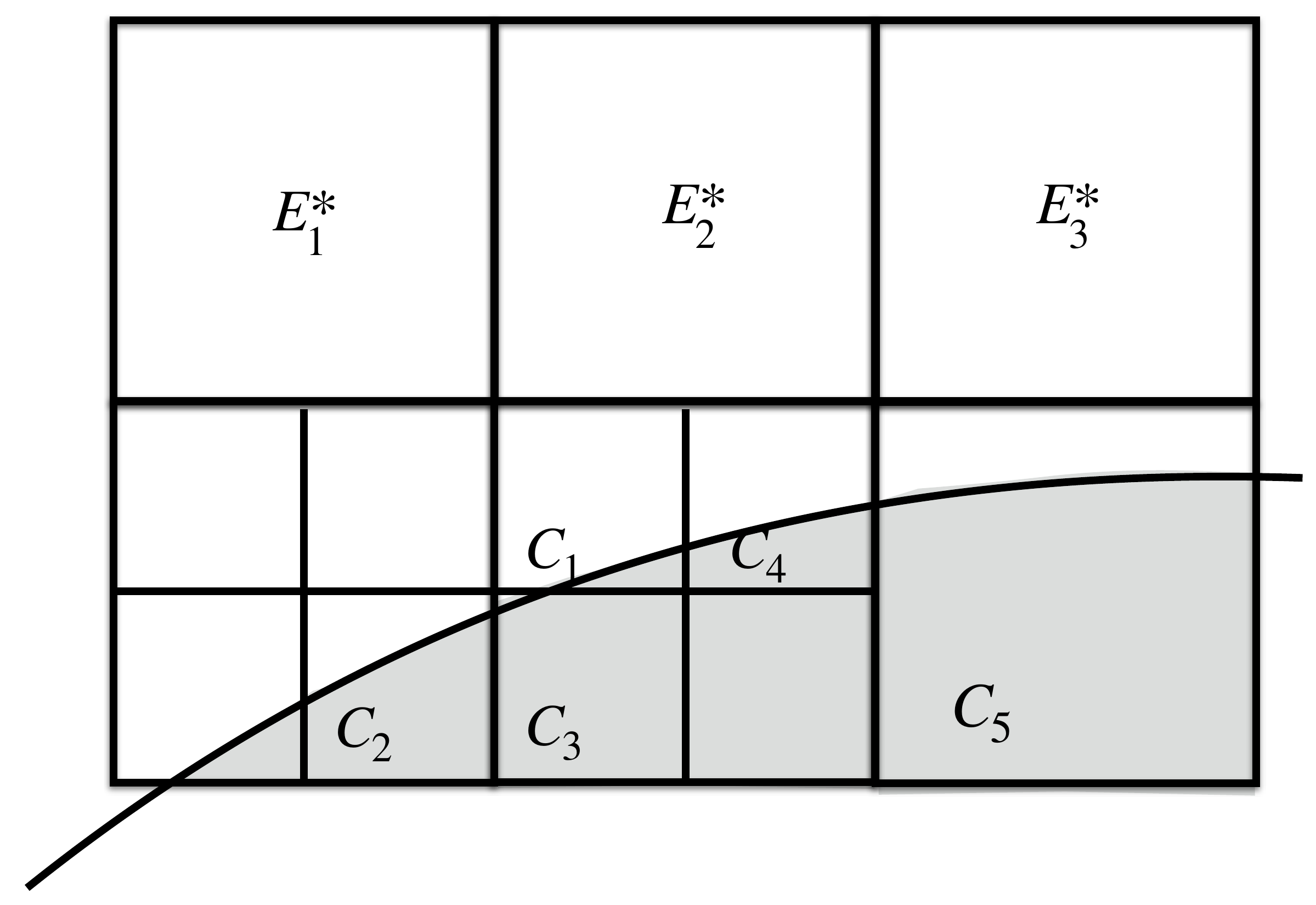} \hspace{.2in}
\includegraphics[width=2in]{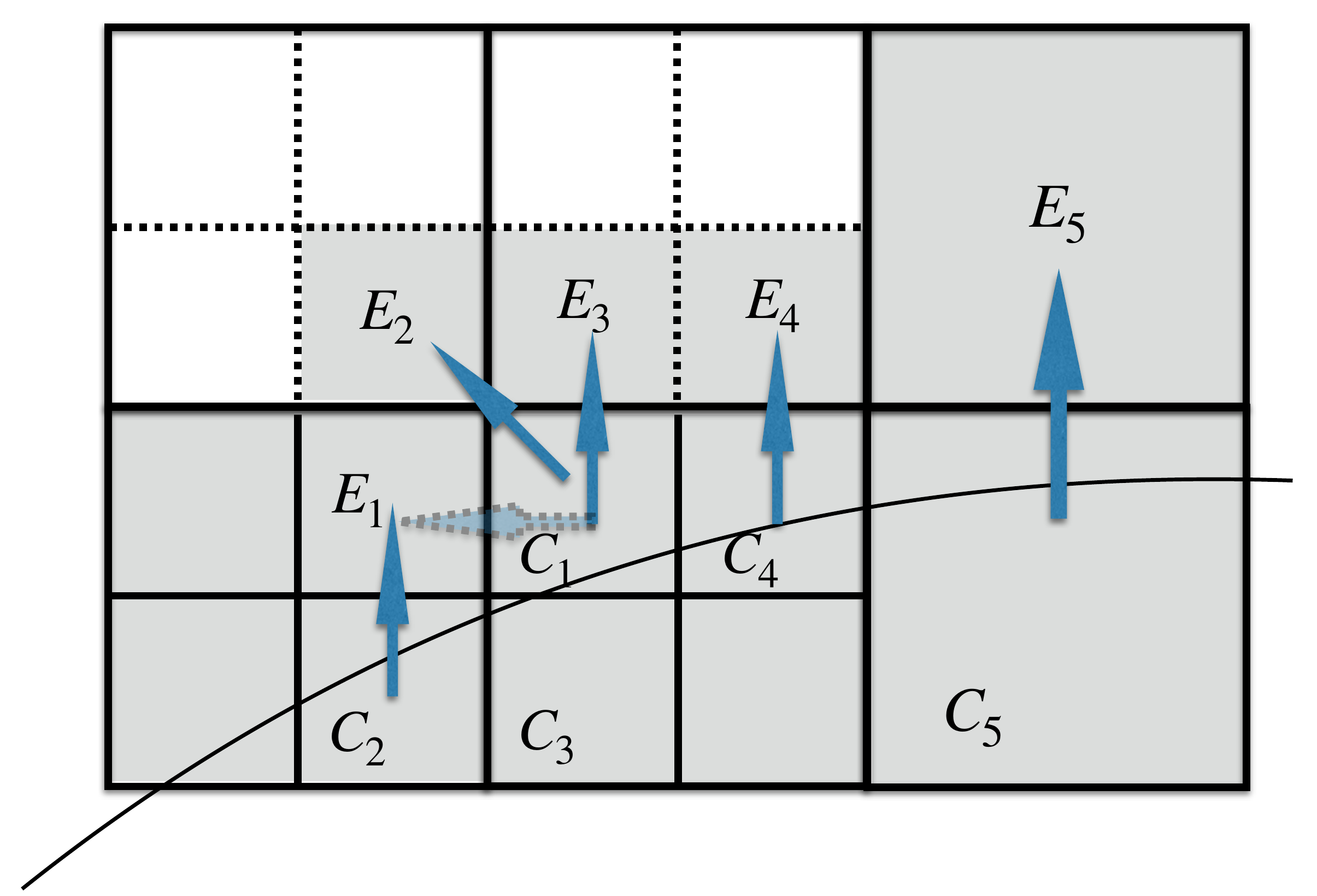} 
 \caption{Cut squares $C_1 - C_5$ on an adaptive,
level-restricted quad-tree and their extensions from
the domain $\Omega$ to a domain $E$ for which the VFMM can be used.
({\em left}): Note that, in the initial discretization, $C_1$ is the
closest box to $E_1^*$ and $E_2^*$, while $C_5$ is closest to 
$E_3^*$. 
({\em right}): Since $C_1$ is at a finer level, the boxes
$E_1^*$ and $E_2^*$ are subdivided before extension, 
while $E_3^*$ is not - it is within the extension region for 
$C_5$. The arrows indicate the box from which the extension to 
the indicated exterior boxes is computed.}
 \label{fig:makext}
\end{figure}

{In addition to ensuring that the source density is resolved, 
we require
two addition conditions to be satisfied on the discretization. First, we 
we impose what we call an  \emph{extension-restriction}, meaning that a
cut leaf square $S$ cannot have coarse neighbors in its extension list
(see \cref{fig:makext}).
Second, For cut leaf squares, we require that the side length
of the box be less than or equal to twice the length of the boundary segment 
$|\Gamma_i|$ which intersects it.}

\begin{remark}
In the VFMM, as in all FMMs, non-neighboring interactions
are approximated in a hierarchical fashion 
using outgoing (multipole) and incoming (local) expansions
with controllable precision. 
Near neighbor interactions, on the other hand, are weakly singular,
and computed using precomputed tables of integrals. The size of these
tables is quite modest because the level-restriction criterion limits
the number of possible configurations that need to be considered and
there are a fixed set of $K^2$ basis functions and target points
that need to be considered on a given leaf node.
We refer the reader to 
\cite{vfmm,greengard1996direct} and the references therein for details.
\end{remark}

\subsection{Function extension on a quad-tree and the VFMM}
\label{ssec:localext}

Suppose now that $\bar{S}_{i}$ is the extension square associated with the 
cut leaf square $S_{i}$. 
Using the method of \cref{sec:leafextend},
we obtain the $(\bar{K}-1)$th total order Chebyshev expansion of 
$f^{e}$ on $\bar{S}_i$ in the form
\begin{equation}
  \label{eq:chebexpbarS}
  f^{e}(\mathbf{x}) = \sum\limits_{m+n<\bar{K}}\alpha^{i}_{m,n}T_{m}(x_{1})T_{n}(x_{2}),\quad \mathbf{x} = (x_{1},x_{2}) \in \bar{S},
\end{equation}
where $T_m(x)$ is the Chebyshev polynomial of degree $m$ 
scaled to the dimensions of 
$\bar{S}_i$. We then evaluate the expression \cref{eq:chebexpbarS}
at every square in the extension list of $S_i$.

We carry out this procedure for all cut leaf squares in the discretization
(\cref{fig:makext}).
For leaf squares that don't intersect the domain $\Omega$ and are 
not in any cut square's extension list, we set $f^{e}$ to zero.
The set of all regular leaf nodes, all cut leaf nodes (to which $f$ has been 
extended) and all extension squares defines the domain $E$ with non-zero
data, to which the 
VFMM from \cite{vfmm} is applied, computing $\mathcal{V}_{E}[f^{e}]$ on the
$K \times K$ Chebyshev grids for all leaf squares. 
At any point $\mathbf{x}$ in the closure of $\Omega$, 
it is straightforward to compute
$\mathcal{V}_{E}[f^{e}](\mathbf{x})$ by interpolation of the Chebyshev
expansion of the leaf node containing the point.

\section{Boundary correction using a double layer potential}
\label{sec:dblpot}

Having found a particular solution to the Poisson equation in $\Omega$,
namely $\mathcal{V}_{E}[f^{e}]$, it remains to solve
the Laplace equation \eqref{eq:laplaceeq} with modified Dirichlet data:
$g(\mathbf{x})- \mathcal{V}_{E}[f^{e}](\mathbf{x})$.
The contribution $g(\mathbf{x})$ is given by the user and we compute 
the contribution $\mathcal{V}_{E}[f^{e}](\mathbf{x})$ as described in the 
preceding section.
For the remainder of this section we consider the interior problem
for a simply connected domain. 
The modifications required to handle multiply connected domain or exterior 
problems are discussed in \cref{sec:mathprel}.

We solve \eqref{eq:laplaceeq} using the
boundary integral equation \cref{eq:inteq} 
with a Nystr\"{o}m discretization \cite{atkinson_1997}. 
For this, 
let $\{t^{G}_{j}\}_{j=1}^{N_{q}}$ and $\{w^{G}_{j}\}_{j=1}^{N_{q}}$ 
be the canonical Gauss-Legendre nodes and weights for the interval $[-1,1]$. 
Consider a panel $\Gamma_i$ in \eqref{bdrydef}, parametrized as in 
\eqref{panelparam}.
We let 
$\mathbf{y}_{ij} = \mathbf{y}(\boldsymbol{\gamma}_{i}(t^{G}_{j}))$, 
$\boldsymbol{\nu}_{ij} = \boldsymbol{\nu}(\mathbf{y}_{ij})$, 
$s_{ij} = \|\boldsymbol{\gamma}_{i}^{\prime}(t^{G}_{j})\|$, and 
$\sigma_{ij}$ be the approximation of $\sigma(\mathbf{y}_{ij})$. 
Applying Gauss-Legendre quadrature to the double layer potential yields
\begin{equation}
\label{eq:discdblpotpanel}
\mathcal{D}[\sigma](\mathbf{x}) = 
\int\limits_{\partial\Omega}\frac{\partial 
G(\mathbf{x}-\mathbf{y})}{\partial \boldsymbol{\nu}(\mathbf{y})}
\sigma(\mathbf{y})\,\mathrm{d}s
\approx 
\sum\limits_{i = 1}^{N_{\Gamma}}\sum\limits_{j = 1}^{N_{q}}\frac{\partial G(\mathbf{x},\mathbf{y}_{ij})}{\partial \boldsymbol{\nu}_{ij}}\sigma_{ij}s_{ij}w^{G}_{j},
\end{equation}
since $\mathrm{ds} = \|\boldsymbol{\gamma}^{\prime}_{i}(t)\|\,\mathrm{d}t$.
Recall that the double layer kernel is smooth on a smooth boundary, so that
the approximation \cref{eq:discdblpotpanel} has an error of the order
$O(h^{2N_{q}}$ where $h = |\Gamma_i|/N_q$.
Using this quadrature in our Nystr\"{o}m scheme applied to \cref{eq:inteq} 
yields the discrete linear system
\begin{equation}
  \label{eq:discinteq}
  \frac{1}{2}\sigma_{i^{\prime}j^{\prime}} + \frac{1}{2\pi}\sum\limits_{i=1}^{N_{\Gamma}}\sum\limits_{j=1}^{N_{q}}\frac{\partial G(\mathbf{y}_{i^{\prime}j^{\prime}},\mathbf{y}_{ij})}{\partial \boldsymbol{\nu}_{ij}}\sigma_{ij}s_{ij}w^{G}_{j} = g(\mathbf{y}_{i^{\prime}j^{\prime}})-\mathcal{V}_{E}[f^{e}](\mathbf{y}_{i^{\prime}j^{\prime}}),
\end{equation}
for $i^{\prime} = 1,\ldots,N_{\Gamma}$ and $j^{\prime} = 1,\ldots,N_{q}$. 
In matrix form, we write \eqref{eq:discinteq} as
\[
(\boldsymbol{\mathrm{I}} + 2 \mathbf{D})\boldsymbol{\sigma} 
= 2(\mathbf{g}-\mathbf{V}).
\] 
While the system matrix $\mathbf{D}$ is dense, it is well-conditioned
and can be solved efficiently with GMRES. This follows from
the fact that the underlying
integral equation is an invertible Fredholm equation of the second kind,
whose eigenvalues cluster at $(1,0)$ 
\cite{2007introductory,KressRainer2014LIE,trefethen1997numerical}.) 
Furthermore, the matrix-vector multiplications required by GMRES can be 
computed using the original (``point") FMM 
with $\mathcal{O}(N_{\Gamma}N_{q})$ operations, resulting in an optimal
time solver \cite{GREENBAUM1993267,greengard1987fast,R-JCP-1985}.

Having solved the integral equation, we may
evaluate the double layer potential \cref{eq:discdblpotpanel} at all
interior points using the point FMM \cite{greengard1987fast}. 
Care must be taken, however, 
as $\mathbf{x}$ approaches the boundary $\partial \Omega$,
since the kernel \cref{eq:dblker} is singular and
the smooth Gauss-Legendre rule used above loses accuracy.
Designing quadrature rules for this regime has been an active area of 
research, and 
there are several FMM-compatible methods available that restore precision, 
such as \cite{KLOCKNER2013332,ABarnettquad,BarnettVeeraWuquad}. 
We use the Helsing-Ojala correction scheme \cite{HelsingOjala} in this paper.

\subsection{Error analysis}

One of the advantages of potential theory is that it uncouples 
the discretization of the domain from that of the boundary and permits
very simple error analysis.
In computing the particular solution $\mathcal{V}_{E}[f^{e}]$, 
there are two sources of error.
The first is the error $\epsilon_f$ made in the piecewise polynomial
approximation of $f(\mathbf{x})$.
Since the volume integral operator $\mathcal{V}_{E}$ 
is bounded, this contributes an error of the order $O(\epsilon_f)$.
The second is the error made in computing $\mathcal{V}_{E}[f^{e}]$ 
for that piecewise polynomial approximation $f^e$.
The VFMM computes this exactly, up to the tolerance $\epsilon_{FMM}$
specified by the user.
More complicated is the error in the double layer potential. 
Since this involves the solution of an integral equation, 
we can't specify the accuracy {\em a priori}. We can say, however, that
the order of accuracy of the solution is that of the underlying quadrature
rule. This is a particular
feature of second kind integral equations \cite{atkinson_1997}.
That is, we are guaranteed high order convergence from a high order 
accurate rule. We must also ensure that 
the right-hand side of the integral equation
\eqref{eq:discinteq}, namely
$g(\mathbf{y})-\mathcal{V}_{E}[f^{e}](\mathbf{y})$,
is well-resolved. This is a slightly
subtle issue, since the function is cut off sharply at the
boundary of the extension region $E$, which could introduce
high-frequency content in the term
$\mathcal{V}_{E}[f^{e}](\mathbf{y})$.
Our algorithm mitigates this by ensuring that 
the corner points of the polygonal boundary
$\partial E$ are pushed out at least a full leaf node away from 
the domain boundary $\partial \Omega$.

\begin{remark}
One could also sample the curve more finely to ensure that a
piecewise polynomial approximation of 
$g(\mathbf{y})-\mathcal{V}_{E}[f^{e}](\mathbf{y})$ is resolved
to the desired precision. We have not investigated this issue in detail.
In the present paper, we sample the boundary
sufficiently finely that the error is dominated by the
accuracy of the volume integral. 
\end{remark}

\begin{remark}
For nonsmooth boundaries, we rely on the recent development of
high order solvers that deal efficiently with corner singularities,
such as \cite{brs,rciphelsing,hrs2019}). 
The essential idea in these schemes 
is the use of dyadic refinement to the corner to overcome the induced 
singularity in the double layer density.
We make use here of the RCIP method 
of \cite{rciphelsing} and refer the reader to the original paper for 
further details. 
\end{remark}

\section{Numerical results}
\label{sec:numres}

The bulk of the software for our function extension scheme 
is written in Julia 1.7.1 \cite{bezanson2017julia} and available 
at \cite{bieps2d}.
It can be used to generate the results in this section. 
Software for the boundary integral equation, 
the evaluation of the double layer potential, the RCIP scheme, 
and the RBF-QR algorithm 
have also been implemented in Julia. The latter is available
at \cite{rbfqrcode}. 
The full Poisson solver relies on several external packages: 
the VFMM \cite{vfmm} is written in Fortran 
and available at \cite{fmmvolgithub}, fixed at eighth order accuracy.
We set the FMM tolerance to $\epsilon_{FMM} = 0.5 \times 10^{-11}$.
The ``point" FMM we use is available at \cite{fmmlib2d}. 

In our discretization, we set $K=8$ for the Chebyshev grids on leaf nodes,
whether they are regular or cut. The number of Gaussians is set to 
$N_{\mathcal{G}}=66$, as discussed in \cref{sec:leafextend}.
We have found this works well in practice to obtain  eighth order
accuracy.
On the extension region, $\bar{S}$, we set 
$\bar{K} = 12$. When it is resampled on the individual extension squares,
however, we interpolate on $8 \times 8$ Chebyshev grids, for 
compatibility with the VFMM.
As the algorithm traverses the extension list, no square is written
to more than once, making the extension step trivially parallel. 
On the boundary, we use $N_{q} = 16$ Gauss-Legendre nodes for each panel. 
The number of panels $N_{\Gamma}$ is set to be sufficiently large that 
resolving the geometry does not dominate the error. 
That is, we pick $N_{\Gamma}$ to ensure that, on each panel, 
the $16$ point Gauss-Legendre expansion of $\|\boldsymbol{\gamma}\|$ is 
resolved to fifteen digits of accuracy.

In the following numerical experiments, we compute the solution at the 
subset of a uniform $100 \times 100$ grid on $D$ that lie inside $\Omega$ 
for the interior problem, and outside $\Omega$ for the exterior problem. 
We measure the error in both the relative $\ell^{\infty}$ norm and the 
relative $\ell^{2}$ norm. 

For convergence studies, we use a uniform quad-tree; thus,
 at level $l$ there are $N =  8 \cdot 2^{l}$ points in each dimension. 
{All computations were carried out on a single core of a 
$4.2$ GHz Intel i$7-8620U$ with $16$ GB of memory.}
 
\subsection{The interior problem}
\label{ssec:intpois}
\begin{figure}[h!]
\centering
    \subfloat[\centering ]{{\includegraphics[trim=0.0cm 0.0cm 1.0cm 0.0cm, clip,width=4.5cm]{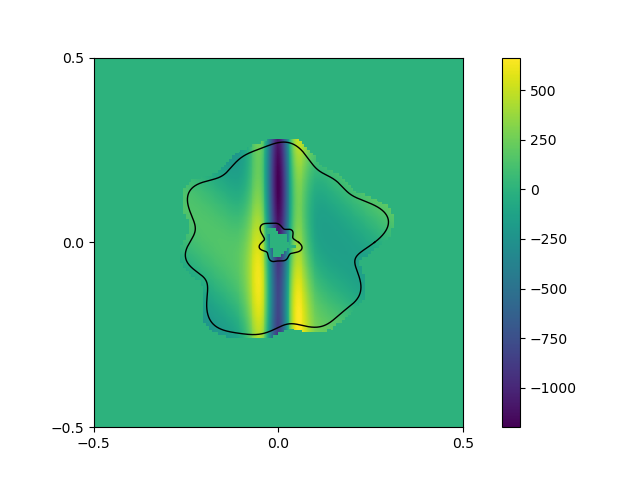} }}%
    \subfloat[\centering ]{{\includegraphics[trim=0.0cm 0.0cm 1.0cm 0.0cm, clip,width=4.5cm]{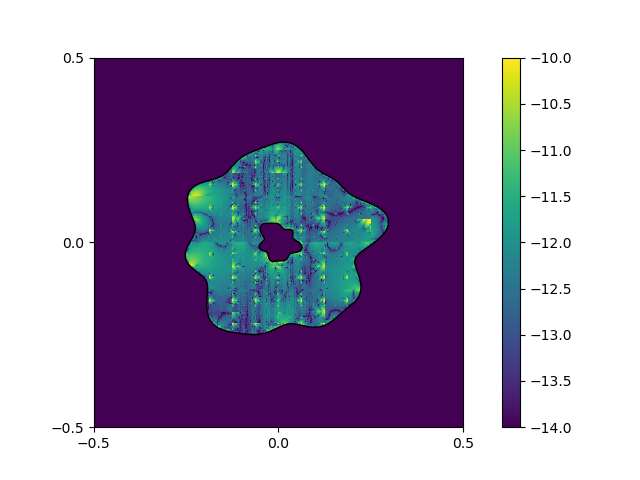}}}%
    \subfloat[\centering ]{{
%
%
\definecolor{mycolor1}{rgb}{0.00000,0.44700,0.74100}%
\definecolor{mycolor2}{rgb}{0.85000,0.32500,0.09800}%
\definecolor{mycolor3}{rgb}{0.92900,0.69400,0.12500}%
\definecolor{mycolor4}{rgb}{0.49400,0.18400,0.55600}%
\definecolor{mycolor5}{rgb}{0.46600,0.67400,0.18800}%
\definecolor{mycolor6}{rgb}{0.00000,0.74902,0.74902}%
\begin{tikzpicture}[scale=0.47]

\begin{loglogaxis}[%
xticklabel style={/pgf/number format/fixed},
xmin=0,
xmax=2500,
xlabel style={font=\color{white!15!black}},
xlabel={$N$},
xtick={50,100,200,400,800,1600},
xticklabels={$50$,$100$,$200$,$400$,$800$,$1600$},
xticklabel style = {font=\small},
xminorticks=false,
ymin=1e-13,
ymax=1,
yminorticks=true,
ylabel style={font=\small\color{white!15!black}},
ylabel={Relative error},
ytick={1,1e-1,1e-2,1e-3,1e-4,1e-5,1e-6,1e-7,1e-8,1e-9,1e-10,1e-11,1e-12,1e-13,1e-14},
yticklabels={$10^{0}$,$$,$10^{-2}$,$$,$10^{-4}$,$$,$10^{-6}$,$$,$10^{-8}$,$$,$10^{-10}$,$$,$10^{-12}$,$$,$10^{-14}$},
yticklabel style = {font=\small},
axis background/.style={fill=white},
xmajorgrids,
xminorgrids,
ymajorgrids,
legend style={font=\normalsize,at={(1.0,1.0)}, anchor=north east, legend cell align=left, align=left, draw=white!15!black},
clip mode=individual,
]

\addplot [color=mycolor4, line width=1.0pt, mark=triangle, mark options={solid, mycolor4}, mark size=2.0pt]
  table[row sep=crcr]{%
64.0 0.042566411535742844\\
128.0 0.0001317230021190127\\
256.0 9.882552819353315e-7\\
512.0 2.4408239247874183e-10\\
1024.0 1.6952608687337426e-11\\
2048.0 1.7020505198538133e-11\\
};
\addlegendentry{$\ell^{\infty}$}

\addplot [color=mycolor3, line width=1.0pt, mark=*, mark options={solid, mycolor3}, mark size=2.0pt]
  table[row sep=crcr]{%
64.0 0.0007772000476552755\\
128.0 5.8135359733672644e-6\\
256.0 2.576837610850237e-9\\
512.0 2.4298232058961563e-11\\
1024.0 1.7030708226177585e-11\\
2048.0 1.704072574422359e-11\\
};
\addlegendentry{Analytic, $\ell^{\infty}$}

\addplot [color=black, line width=1.0pt, mark options={solid, black}]
  table[row sep=crcr]{%
64.0 0.004472601531070688\\
128.0 1.7471099730744875e-5\\
256.0 6.824648332322217e-8\\
512.0 2.665878254813366e-10\\
1024.0 1.041358693286471e-12\\
2048.0 4.0678073956502776e-15\\
};
\addlegendentry{$8$th order}

\addplot [color=mycolor2, line width=1.0pt, mark=+, mark options={solid, mycolor2}, mark size=3.0pt]
  table[row sep=crcr]{%
64.0 0.002675061508370936\\
128.0 6.776708317503218e-6\\
256.0 5.17741104087048e-8\\
512.0 1.42759849174401e-11\\
1024.0 8.540780691755239e-13\\
2048.0 8.048971585769807e-13\\
};
\addlegendentry{$\ell^{2}$}

\addplot [color=mycolor1, line width=1.0pt, mark=square*, mark options={solid, mycolor1}, mark size=2.0pt]
  table[row sep=crcr]{%
64.0 9.253483436714584e-5\\
128.0 7.667288572692666e-7\\
256.0 8.406571122114159e-10\\
512.0 8.232719688338623e-12\\
1024.0 8.225801431611441e-13\\
2048.0 8.039454215999703e-13\\
};
\addlegendentry{Analytic, $\ell^{2}$}

\end{loglogaxis}
\end{tikzpicture}%

}}%
    \caption{
The doubly-connected interior problem of \cref{ssec:intpois}.
(a) The extended function $f = \Delta u$ 
with $u$ given by \cref{eq:harmext} (with seven levels of refinement). 
(b) Pointwise error in the computed solution. 
(c) Convergence plot under uniform refinement (markers
are from levels three to eight in the refinement process).}
    \label{fig:interior}%
\end{figure}
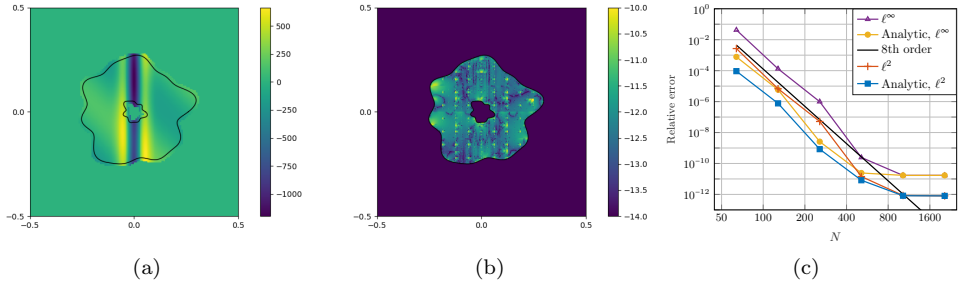

For our first test, we consider the problem posed as
Example $2$ in \cite{ASKHAM20171}. 
It involves a doubly connected domain with a right-hand side that has 
some very fine features with exact solution 
\begin{equation}
  \label{eq:harmext}
  u(\mathbf{x}) = \sin(10(x_{1} + x_{2})) + x_{1}^2 - 3x_{2} + 8 + \exp(-500x_{1}^2), \quad \mathbf{x} \in \Omega.
\end{equation}
The two boundary components are specified in polar coordinates 
with $\theta \in [0,2\pi)$ and 
\[ r(\theta) = \sum_{j}(c_{j}\cos(j\theta) + d_{j}\sin(j\theta)).
\] 
The non-zero coefficients for the outer boundary $\partial \Omega_{0}$ are $c_{0} = 0.25$, $d_{3} = c_{6} = c_{8} = c_{10} = 0.01$ and $c_{5} = 0.02$. 
The non-zero coefficients for the inner boundary $\partial \Omega_{1}$ are $c_{0} = 0.05$ and $c_{2} = d_{3} = c_{5} = c_{7} = 0.005$. 
(See \cref{fig:interior}.) We discretize $\partial \Omega_{0}$ 
with $200$ panels and $\partial \Omega_{1}$ with $180$ panels.

In \cref{fig:interior}, we observe the expected eighth order convergence 
as we refine the quad-tree  uniformly. The $\ell^{\infty}$ error levels 
out after six levels of refinement at about eleven digits of accuracy, more or less the FMM tolerance $\epsilon_{FMM}$. The $\ell^{2}$ error continues to
decrease for one more level, reaching twelve digits of accuracy.
For comparison, we also plot the 
errors when the function extension is carried out exactly based on the 
exact solution (to the same region $E \supset \Omega$). 
We refer to this as the \emph{analytic extension}. Note that 
we lose one to two digits of accuracy from our numerical scheme
(although with sufficient refinement, the errors are the same).

We test the performance of the adaptive solver, using 
$\epsilon_{FMM} = 0.5\times 10^{-11}$ for both the VFMM and in 
determining when the right-hand side is sufficiently resolved.
As noted above, we also ensure that the dimensions of the cut squares are
commensurate with the boundary panel size 
($| \partial\Omega_0|/N_{\Gamma}$), 
which requires seven levels of refinement near the boundary. 
The resulting $\ell^{\infty}$ error is $2\times 10^{-11}$, with 
an $\ell^{2}$ error of $10^{-12}$. The full discretization requires
$3361$ leaf squares, of which $339$ are cut, with a total of about 
$215,000$ points. 
The construction of the quad-tree, which includes labeling squares as cut, 
imposing the level-restriction, and imposing the extension-restriction, 
requires $0.3$ seconds. The precomputation steps in function extension -
building the extension lists and identifying points as inside or outside -
requires $0.7$ seconds. 
Creating the extension itself requires $0.4$ 
seconds, and the VFMM requires $0.2$ seconds. 
Finally, solving the integral equation and evaluating the double layer
potential requires $0.7$ seconds. 
Note that,
since seven levels of uniform refinement would require $N\sim10^{3}$ points, 
adaptivity has yielded a factor of five improvement for the same accuracy.
It is difficult to make a direct comparison with the scheme of 
\cite{ASKHAM20171}, since they used a less smooth extension and relied
on a fourth order VFMM. For the same example, however,
ten times more points were needed
to obtain an error of $10^{-8}$.

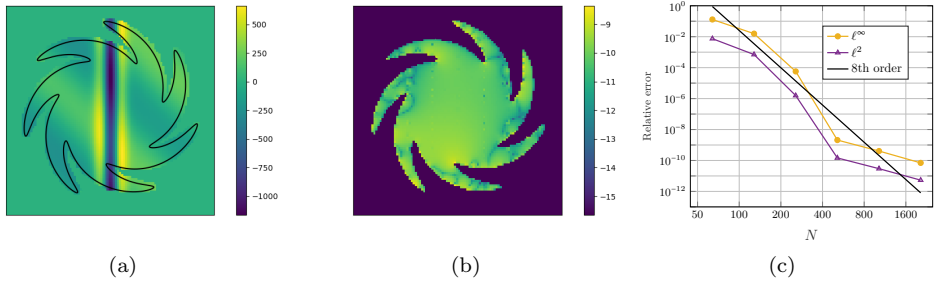
\begin{figure}
\centering
    \subfloat[\centering ]{{\includegraphics[trim=0.0cm 0.0cm 1.0cm 0.0cm, clip,width=4.5cm]{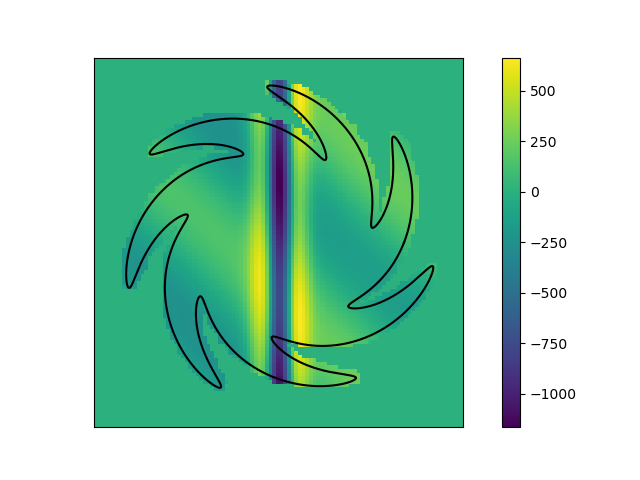} }}%
    \subfloat[\centering ]{{\includegraphics[trim=0.0cm 0.0cm 1.0cm 0.0cm, clip,width=4.5cm]{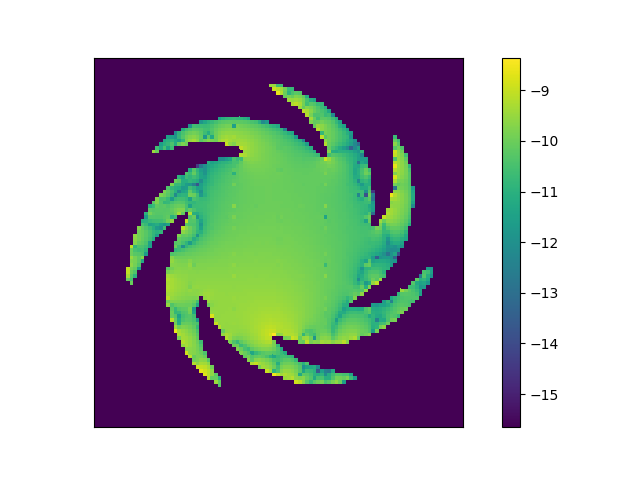}}}%
    \subfloat[\centering ]{{
%
%
\definecolor{mycolor1}{rgb}{0.00000,0.44700,0.74100}%
\definecolor{mycolor2}{rgb}{0.85000,0.32500,0.09800}%
\definecolor{mycolor3}{rgb}{0.92900,0.69400,0.12500}%
\definecolor{mycolor4}{rgb}{0.49400,0.18400,0.55600}%
\definecolor{mycolor5}{rgb}{0.46600,0.67400,0.18800}%
\definecolor{mycolor6}{rgb}{0.00000,0.74902,0.74902}%
\begin{tikzpicture}[scale=0.47]

\begin{loglogaxis}[%
xticklabel style={/pgf/number format/fixed},
xmin=0,
xmax=2500,
xlabel style={font=\large\color{white!15!black}},
xlabel={$N$},
xtick={50,100,200,400,800,1600},
xticklabels={$50$,$100$,$200$,$400$,$800$,$1600$},
xticklabel style = {font=\small},
xminorticks=false,
ymin=1e-13,
ymax=1,
yminorticks=true,
ylabel style={font=\small\color{white!15!black}},
ylabel={Relative error},
ytick={1,1e-1,1e-2,1e-3,1e-4,1e-5,1e-6,1e-7,1e-8,1e-9,1e-10,1e-11,1e-12,1e-13,1e-14},
yticklabels={$10^{0}$,$$,$10^{-2}$,$$,$10^{-4}$,$$,$10^{-6}$,$$,$10^{-8}$,$$,$10^{-10}$,$$,$10^{-12}$,$$,$10^{-14}$},
yticklabel style = {font=\small},
axis background/.style={fill=white},
xmajorgrids,
xminorgrids,
ymajorgrids,
legend style={font=\normalsize,at={(0.9,0.9)}, anchor=north east, legend cell align=left, align=left, draw=white!15!black},
clip mode=individual,
]
\addplot [color=mycolor3, line width=1.0pt, mark=*, mark options={solid, mycolor3}, mark size=2.0pt]
  table[row sep=crcr]{%
64.0 0.12974342782621415\\
128.0 0.015737270620860867\\
256.0 5.710327509581139e-5\\
512.0 2.0870233892640344e-9\\
1024.0 4.0780997723454013e-10\\
2048.0 6.980751138622477e-11\\
};
\addlegendentry{$\ell^{\infty}$}

\addplot [color=mycolor4, line width=1.0pt, mark=triangle, mark options={solid, mycolor4}, mark size=2.0pt]
  table[row sep=crcr]{%
64.0 0.0074551442960338775\\
128.0 0.0007023735635492773\\
256.0 1.5654520669559991e-6\\
512.0 1.458419619978537e-10\\
1024.0 2.9538904099482895e-11\\
2048.0 5.3507490432581515e-12\\
};
\addlegendentry{$\ell^{2}$}

\addplot [color=black, line width=1.0pt, mark options={solid, black}]
  table[row sep=crcr]{%
64.0 0.892401328481747\\
128.0 0.0034859426893818243\\
256.0 1.3616963630397751e-5\\
512.0 5.3191264181241216e-8\\
1024.0 2.077783757079735e-10\\
2048.0 8.116342801092715e-13\\
};
\addlegendentry{$8$th order}

\end{loglogaxis}
\end{tikzpicture}%


}}%
    \caption{(a) Extended function $f = \Delta u$, with $u$ from \cref{eq:harmext}, using seven levels of refinement. (b) Pointwise absolute error for solving the problem in \cref{ssec:intpois}. (c) Convergence plot as the quad-tree is uniformly refined. Here, we plot from three levels of refinement to eight levels of refinement.}%
    \label{fig:interiorsaw}%
\end{figure}

In a second test, we use the same exact solution
$u$ from \cref{eq:harmext}, but in the simply connected domain shown 
in \cref{fig:interiorsaw}. Using complex notation 
$z(\theta) = x(\theta) + i y(\theta)$, the boundary $\partial \Omega$
is given by
\begin{equation}
  z = 0.17((2 + 0.5 \sin(7\theta))\cos(\theta + 0.5\sin(7\theta)) + i((2 + 0.5\sin(7\theta))\sin(\theta + 0.5\sin(7\theta))),
\end{equation}
where $\theta \in [0,2\pi)$. We again observe the expected eighth order 
convergence, but with a larger constant for the error than for our first example.

\subsection{The exterior problem}
\label{ssec:extpoiss}

\begin{figure}
\centering
    \subfloat[\centering ]{{\includegraphics[trim=0.0cm 0.0cm 1.0cm 0.0cm, clip,width=4.5cm]{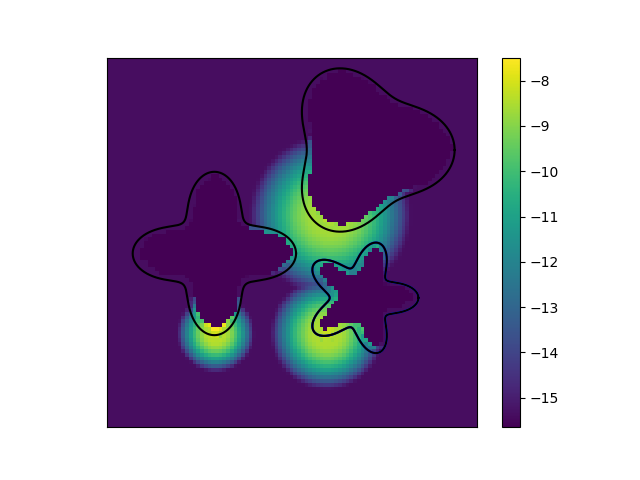} }}%
    \subfloat[\centering ]{{\includegraphics[trim=0.0cm 0.0cm 1.0cm 0.0cm, clip,width=4.5cm]{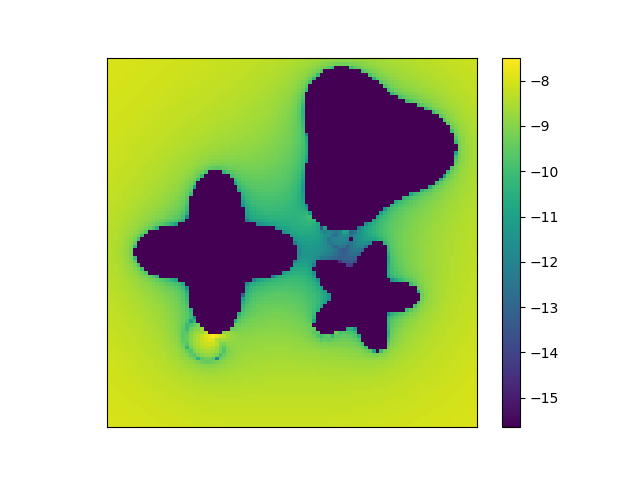}}}%
    \subfloat[\centering ]{{
%
%
\definecolor{mycolor1}{rgb}{0.00000,0.44700,0.74100}%
\definecolor{mycolor2}{rgb}{0.85000,0.32500,0.09800}%
\definecolor{mycolor3}{rgb}{0.92900,0.69400,0.12500}%
\definecolor{mycolor4}{rgb}{0.49400,0.18400,0.55600}%
\definecolor{mycolor5}{rgb}{0.46600,0.67400,0.18800}%
\definecolor{mycolor6}{rgb}{0.00000,0.74902,0.74902}%
\begin{tikzpicture}[scale=0.468]

\begin{loglogaxis}[%
xticklabel style={/pgf/number format/fixed},
xmin=0,
xmax=4200,
xlabel style={font=\small\color{white!15!black}},
xlabel={$N$},
xtick={50,100,200,400,800,1600,3200},
xticklabels={$50$,$100$,$200$,$400$,$800$,$1600$,$3200$},
xticklabel style = {font=\small},
xminorticks=false,
ymin=1e-13,
ymax=1,
yminorticks=true,
ylabel style={font=\small\color{white!15!black}},
ylabel={Relative error},
ytick={1,1e-1,1e-2,1e-3,1e-4,1e-5,1e-6,1e-7,1e-8,1e-9,1e-10,1e-11,1e-12,1e-13,1e-14},
yticklabels={$10^{0}$,$$,$10^{-2}$,$$,$10^{-4}$,$$,$10^{-6}$,$$,$10^{-8}$,$$,$10^{-10}$,$$,$10^{-12}$,$$,$10^{-14}$},
yticklabel style = {font=\small},
axis background/.style={fill=white},
xmajorgrids,
xminorgrids,
ymajorgrids,
legend style={font=\normalsize,at={(0.96,0.96)}, anchor=north east, legend cell align=left, align=left, draw=white!15!black},
clip mode=individual,
]
\addplot [color=mycolor3, line width=1.0pt, mark=*, mark options={solid, mycolor3}, mark size=2.0pt]
  table[row sep=crcr]{%
64.0 0.0027195579561574566\\
128.0 0.067622529565957\\
256.0 0.0002594783513016564\\
512.0 4.857923236627321e-6\\
1024.0 6.795486996283798e-10\\
2048.0 6.5103439663172905e-12\\
4096.0 1.1117621869940231e-12\\
};
\addlegendentry{$\ell^{\infty}$}

\addplot [color=mycolor4, line width=1.0pt, mark=triangle, mark options={solid, mycolor4}, mark size=2.0pt]
  table[row sep=crcr]{%
64.0 0.0004077840210782966\\
128.0 0.003805297489714605\\
256.0 5.899349818384496e-5\\
512.0 4.887183642208493e-6\\
1024.0 2.3835803842159246e-10\\
2048.0 6.53345569709118e-12\\
4096.0 1.1172000138538937e-12\\
};
\addlegendentry{$\ell^{2}$}

\addplot [color=black, line width=1.0pt, mark options={solid, black}]
  table[row sep=crcr]{%
64.0 8.92401328481747\\
128.0 0.03485942689381824\\
256.0 0.0001361696363039775\\
512.0 5.319126418124121e-7\\
1024.0 2.077783757079735e-9\\
2048.0 8.116342801092714e-12\\
};
\addlegendentry{$8$th order}

\end{loglogaxis}
\end{tikzpicture}%

}}%
    \caption{(a) Extended function $f = \Delta u$ on a logarithmic scale, 
with $u$ from \cref{eq:exterior}, using seven levels of refinement. 
(b) Pointwise absolute error for solving the problem in \cref{ssec:extpoiss}. 
(c) Convergence plot as the quad-tree is uniformly refined, with 
markers at refinement levels three to nine.}%
    \label{fig:exterior}%
\end{figure}
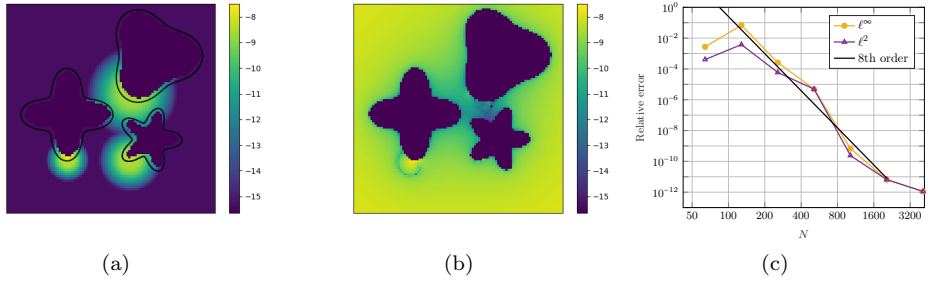

We turn now to a test for the exterior solver on a multiply connected domain 
(\cref{fig:exterior}). Using complex notation again, we have
\begin{equation}
  \label{eq:paramstarfish}
  z = R((1 + a\cos(Nt))\exp(-it)) + c_{1} + ic_{2}, \quad \text{ for } t\in[0,2\pi). 
\end{equation}
For $\partial\Omega_{1}$ we set $R = 0.12$, $a = 0.3$, $N = 5$, $c_{1} = 0.186$, and $c_{2} = -0.15$. For $\partial\Omega_{2}$ we set 
$R = 0.17$, $a = 0.3$, $N = 4$, $c_{1} = -0.21$, and $c_{2} = -0.03$, 
and for $\partial\Omega_{3}$ we set $R = 0.2$, $a = 0.2$, $N = 3$, $c_{1}= 0.2$, and $c_{2} = 0.15$. 
We solve the Poisson equation on this domain with the exact solution
\begin{equation}
  \label{eq:exterior}
  u(\mathbf{x}) = \sum\limits_{j = 1}^{3}e^{-\|\mathbf{x} - \mathbf{y}_{j}\|^{2} / \beta_{j}} -10\log\left(\sqrt{(x_{1} + 0.2)^{2} + x_{2}^{2}}\right), \quad \text{ for } \mathbf{x} \in \mathbb{R}^{2}\setminus(\bar{\Omega}_{1}\cup\bar{\Omega}_{2}\cup\bar{\Omega}_{3}),
\end{equation}
where $\mathbf{y}_{1} = (0.1,0.07)$, $\beta_{1}=10^{-3}$, $\mathbf{y}_{2} = (0.09,-0.25)$, $\beta_{2}=10^{-3}/2.1$, $\mathbf{y}_{3} = (-0.21,-0.25)$, and 
$\beta_{3}=10^{-3}/4.5$ (see \cref{fig:exterior}.) 
Note that the Gaussian centers $\{\mathbf{y}_{j}\}_{j = 1}^{3}$ are 
interior to but close to the boundaries of the inclusions $\Omega_j$. 
Note also that we are seeking a solution which is growing
logarithmically using the representation 
\eqref{eq:solutionmulte} for our integral equation solver, to impose
the radiation  condition $ u(\mathbf{x}) 
\rightarrow 10 \log \| \mathbf{x} \|$
as $\| \mathbf{x} \| \rightarrow \infty$.
Note, however, that the source distribution $f$ may itself have
net ``charge" 
$A^{e} = \iint_{E}f^{e}(\mathbf{y})\,\mathrm{d}\mathbf{y}$, so that the VFMM
is computing a particular solution with growth 
$A^{e} \log \| \mathbf{x} \|$.
Thus, in our integral equation solver, 
for the constraint conditions \eqref{eq:densconstext}, we enforce
\[
\int\limits_{\partial\Omega_{k}} \sigma(\mathbf{y})\,\mathrm{d}\mathbf{y} = 0,\quad\, k = 1,\ldots,N_{\Omega}-1 \quad , \quad
\sum_{k=1}^{N_{\Omega}} A_{k}= 10 - A^{e}.
\]
The convergence plots in \cref{fig:exterior} show the expected
eighth order convergence under uniform refinement.

\subsection{Piecewise smooth boundaries}
\label{ssec:squarepoiss}

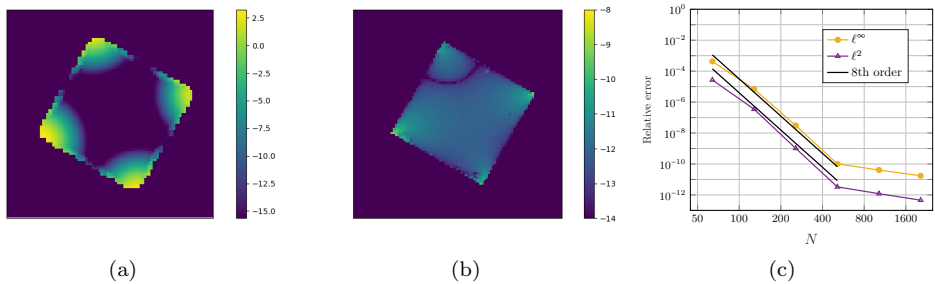
\begin{figure}
\centering
    \subfloat[\centering ]{{\includegraphics[trim=0.0cm 0.0cm 1.0cm 0.0cm, clip,width=4.5cm]{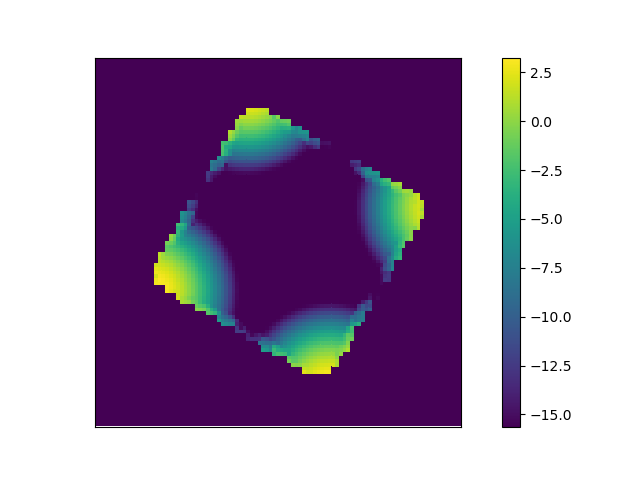} }}%
    \subfloat[\centering ]{{\includegraphics[trim=0.0cm 0.0cm 1.0cm 0.0cm, clip,width=4.5cm]{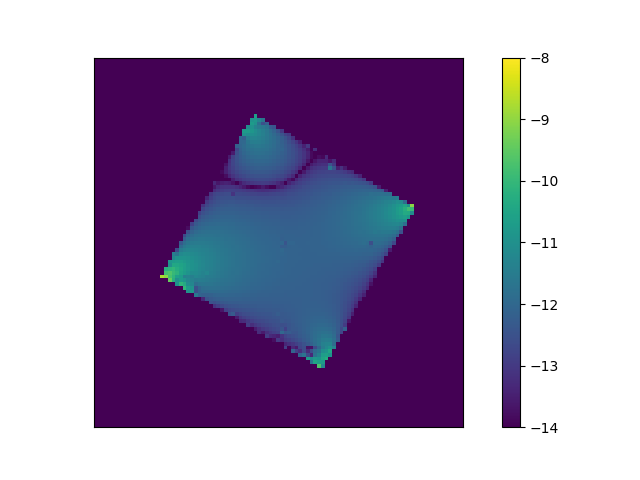}}}%
    \subfloat[\centering ]{{
%
%
\definecolor{mycolor1}{rgb}{0.00000,0.44700,0.74100}%
\definecolor{mycolor2}{rgb}{0.85000,0.32500,0.09800}%
\definecolor{mycolor3}{rgb}{0.92900,0.69400,0.12500}%
\definecolor{mycolor4}{rgb}{0.49400,0.18400,0.55600}%
\definecolor{mycolor5}{rgb}{0.46600,0.67400,0.18800}%
\definecolor{mycolor6}{rgb}{0.00000,0.74902,0.74902}%
\begin{tikzpicture}[scale=0.47]

\begin{loglogaxis}[%
xticklabel style={/pgf/number format/fixed},
xmin=0,
xmax=2500,
xlabel style={font=\large\color{white!15!black}},
xlabel={$N$},
xtick={50,100,200,400,800,1600},
xticklabels={$50$,$100$,$200$,$400$,$800$,$1600$},
xticklabel style = {font=\small},
xminorticks=false,
ymin=1e-13,
ymax=1,
yminorticks=true,
ylabel style={font=\small\color{white!15!black}},
ylabel={Relative error},
ytick={1,1e-1,1e-2,1e-3,1e-4,1e-5,1e-6,1e-7,1e-8,1e-9,1e-10,1e-11,1e-12,1e-13,1e-14},
yticklabels={$10^{0}$,$$,$10^{-2}$,$$,$10^{-4}$,$$,$10^{-6}$,$$,$10^{-8}$,$$,$10^{-10}$,$$,$10^{-12}$,$$,$10^{-14}$},
yticklabel style = {font=\small},
axis background/.style={fill=white},
xmajorgrids,
xminorgrids,
ymajorgrids,
legend style={font=\normalsize,at={(0.9,0.9)}, anchor=north east, legend cell align=left, align=left, draw=white!15!black},
clip mode=individual,
]
\addplot [color=mycolor3, line width=1.0pt, mark=*, mark options={solid, mycolor3}, mark size=2.0pt]
  table[row sep=crcr]{%
64.0 0.00042572151694686455\\
128.0 6.925651564674921e-6\\
256.0 3.055248434088322e-8\\
512.0 1.0192472816553331e-10\\
1024.0 4.0620624406500925e-11\\
2048.0 1.725381765670124e-11\\
};
\addlegendentry{$\ell^{\infty}$}

\addplot [color=mycolor4, line width=1.0pt, mark=triangle, mark options={solid, mycolor4}, mark size=2.0pt]
  table[row sep=crcr]{%
64.0 2.6921612554722883e-5\\
128.0 3.469763407751256e-7\\
256.0 1.0279259891693919e-9\\
512.0 3.3142720949907794e-12\\
1024.0 1.2036569371956432e-12\\
2048.0 4.563939620244778e-13\\
};
\addlegendentry{$\ell^{2}$}

\addplot [color=black, line width=1.0pt, mark options={solid, black}]
  table[row sep=crcr]{%
64.0 0.0011234667099445444\\
128.0 4.3885418357208765e-6\\
256.0 1.7142741545784674e-8\\
512.0 6.696383416322138e-11\\
};
\addlegendentry{$8$th order}

\addplot [color=black, line width=1.0pt, mark options={solid, black}]
  table[row sep=crcr]{%
64.0 0.00014143607904539725\\
128.0 5.52484683771083e-7\\
256.0 2.158143295980793e-9\\
512.0 8.430247249924973e-12\\
};

\end{loglogaxis}
\end{tikzpicture}%

}}%
    \caption{(a) Extended function $f = \Delta u$ on a logarithmic scale, 
with $u$ from \cref{eq:pwsmooth}, using seven levels of refinement. 
(b) Pointwise absolute error for solving the problem in \cref{ssec:extpoiss}. 
(c) Convergence plot as the quad-tree is uniformly refined, with 
markers at refinement levels three to eight.}%
    \label{fig:pwswooth}%
\end{figure}

A good demonstration of the value of potential theory is the solution
of the Poisson equation with a non-smooth boundary. Assuming that the 
source distribution is well-resolved by the user-provided grid, our 
extension scheme is agnostic as to the regularity of the boundary. 
Thus, let us suppose for simplicity that the solution and source density 
are both smooth on a square with side length $0.5$, 
centered at $(0.01,-0.02)$ and rotated $\pi/3$ radians, to avoid
any benefit from alignment with the coordinate axes.
We solve the interior problem with solution
\begin{equation}
  \label{eq:pwsmooth}
  u(\mathbf{x}) = -2\sum\limits_{j = 1}^{4}\left(\mathrm{Ei}(\beta\|\mathbf{x} - \mathbf{x}_{j}\|^{2}) + \log\|\mathbf{x} - \mathbf{x}_{j}\|^{2}\right),
\end{equation}
with  $\beta = 800$, where $\mathrm{Ei}$ denotes the exponential integral 
function, $\mathbf{x}_{1} = (-0.35,-0.135)$, $\mathbf{x}_{2} = (-0.09,042)$, $\mathbf{x}_{3} = (0.445,0.09)$, and $\mathbf{x}_{4} = (0.135,-0.405)$ 
(see \cref{fig:pwswooth}).

No modifications te the code is required for computing the volume potential 
$\mathcal{V}_{D}[f^{e}]$, but the double layer potential develops
a singularity at the corners,so that we require a specialized quadrature
scheme to achieve high order accuracy in computing  the doulbe layer
$\mathcal{D}[\sigma]$. For this, we make use of \emph{recursive(ly) compressed inverse preconditioning} (RCIP) \cite{rciphelsing}. The results are shown 
in \cref{fig:pwswooth}, where we again obtain the expected eighth order convergence. The code works equally well when the solution has corner singularities, 
so long as the source distribution is resolved by the quad-tree.

\subsection{Extension along lines}
\label{ssec:extradial}
\begin{figure}
\centering
    \subfloat[\centering ]{{\includegraphics[trim=10.0cm 2.0cm 1.0cm 0.0cm, clip,width=4.5cm]{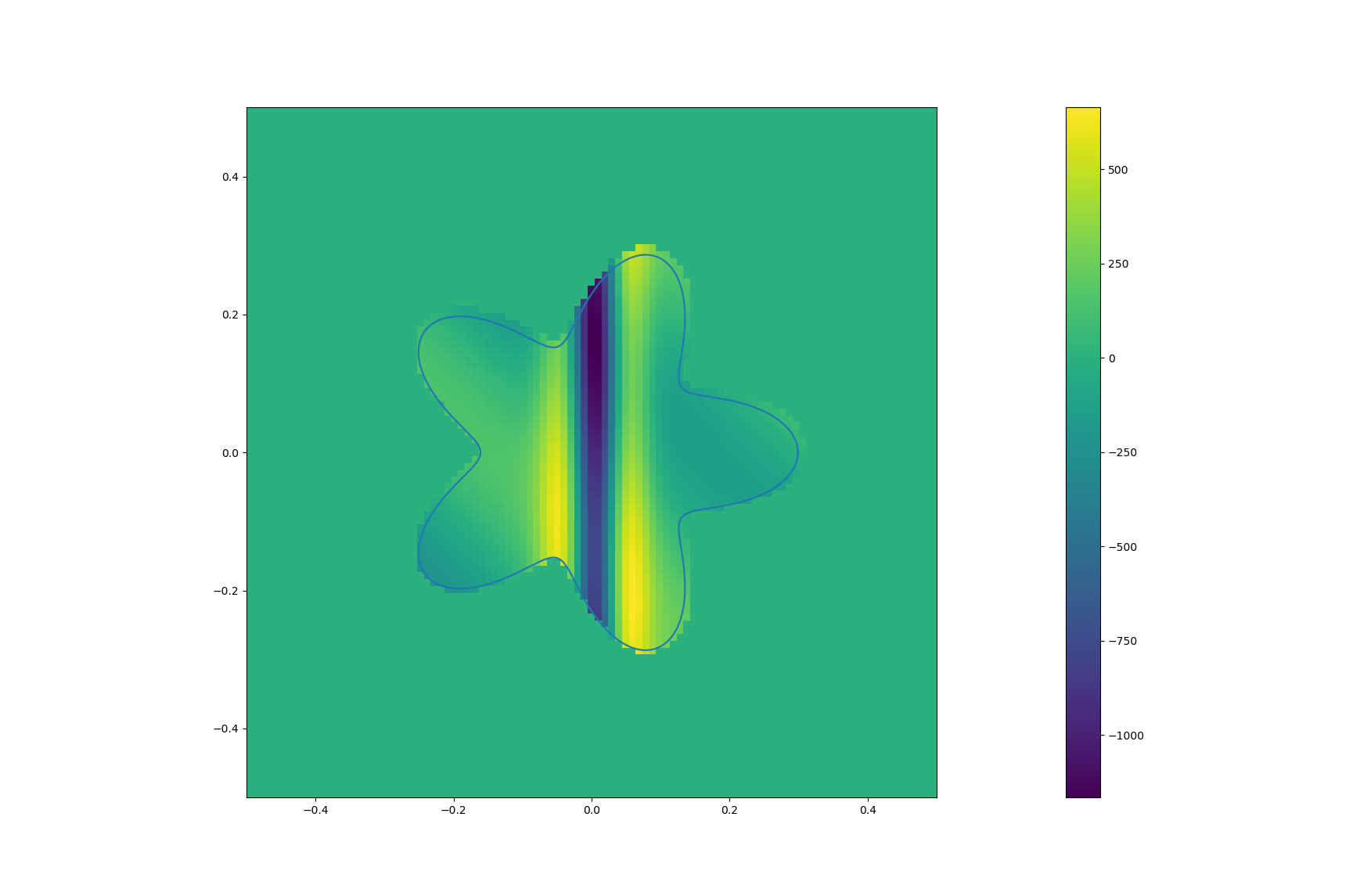} }}%
    \subfloat[\centering ]{{
%
%
\definecolor{mycolor1}{rgb}{0.00000,0.44700,0.74100}%
\definecolor{mycolor2}{rgb}{0.85000,0.32500,0.09800}%
\definecolor{mycolor3}{rgb}{0.92900,0.69400,0.12500}%
\definecolor{mycolor4}{rgb}{0.49400,0.18400,0.55600}%
\definecolor{mycolor5}{rgb}{0.46600,0.67400,0.18800}%
\definecolor{mycolor6}{rgb}{0.00000,0.74902,0.74902}%
\begin{tikzpicture}[scale=0.47]

\begin{loglogaxis}[%
xticklabel style={/pgf/number format/fixed},
xmin=0,
xmax=2500,
xlabel style={font=\small\color{white!15!black}},
xlabel={$N$},
xtick={50,100,200,400,800,1600,2500},
xticklabels={$50$,$100$,$200$,$400$,$800$,$1600$,$2500$},
xticklabel style = {font=\small},
xminorticks=false,
ymin=1e-13,
ymax=1,
yminorticks=true,
ylabel style={font=\small\color{white!15!black}},
ylabel={Relative $\ell^{\infty}$ error},
ytick={1,1e-1,1e-2,1e-3,1e-4,1e-5,1e-6,1e-7,1e-8,1e-9,1e-10,1e-11,1e-12,1e-13,1e-14},
yticklabels={$10^{0}$,$$,$10^{-2}$,$$,$10^{-4}$,$$,$10^{-6}$,$$,$10^{-8}$,$$,$10^{-10}$,$$,$10^{-12}$,$$,$10^{-14}$},
yticklabel style = {font=\small},
axis background/.style={fill=white},
xmajorgrids,
xminorgrids,
ymajorgrids,
legend style={font=\normalsize,at={(0.9,0.9)}, anchor=north east, legend cell align=left, align=left, draw=white!15!black},
clip mode=individual,
]
\addplot [color=mycolor1, line width=1.0pt, mark=o, mark options={solid, mycolor1}, mark size=2.0pt]
  table[row sep=crcr]{%
64.0 0.006636196039715684\\
128.0 0.0004507064057066596\\
256.0 6.250166905148617e-7\\
512.0 7.319212031533491e-11\\
1024.0 1.3053739045115824e-11\\
2048.0 1.3005389015997923e-11\\
};
\addlegendentry{RBF}

\addplot [color=mycolor2, line width=1.0pt, mark=+, mark options={solid, mycolor2}, mark size=2.0pt]
  table[row sep=crcr]{%
64.0 0.021423306587009368\\
128.0 5.808670378743169e-5\\
256.0 2.0655241098596666e-8\\
512.0 1.4445451066993948e-10\\
1024.0 1.3053055652478116e-11\\
2048.0 1.3004791047439928e-11\\
};
\addlegendentry{Rational}


\addplot [color=mycolor4, line width=1.0pt, mark=triangle, mark options={solid, mycolor4}, mark size=2.0pt]
  table[row sep=crcr]{%
64.0 0.0009155321262356562\\
128.0 5.8132763466786175e-5\\
256.0 2.0628164313473924e-8\\
512.0 1.4442973768682255e-10\\
1024.0 1.3053824469195537e-11\\
2048.0 1.3004705623360215e-11\\
};
\addlegendentry{Analytic}

\addplot [color=black, line width=1.0pt, mark options={solid, black}]
  table[row sep=crcr]{%
64.0 0.892401328481747\\
128.0 0.0034859426893818243\\
256.0 1.3616963630397751e-5\\
512.0 5.3191264181241216e-8\\
1024.0 2.077783757079735e-10\\
2048.0 8.116342801092715e-13\\
};
\addlegendentry{$8$th}

\end{loglogaxis}
\end{tikzpicture}%

}}%
    \subfloat[\centering ]{{
%
%
\definecolor{mycolor1}{rgb}{0.00000,0.44700,0.74100}%
\definecolor{mycolor2}{rgb}{0.85000,0.32500,0.09800}%
\definecolor{mycolor3}{rgb}{0.92900,0.69400,0.12500}%
\definecolor{mycolor4}{rgb}{0.49400,0.18400,0.55600}%
\definecolor{mycolor5}{rgb}{0.46600,0.67400,0.18800}%
\definecolor{mycolor6}{rgb}{0.00000,0.74902,0.74902}%
\begin{tikzpicture}[scale=0.47]

\begin{loglogaxis}[%
xticklabel style={/pgf/number format/fixed},
xmin=0,
xmax=2500,
xlabel style={font=\small\color{white!15!black}},
xlabel={$N$},
xtick={50,100,200,400,800,1600,2500},
xticklabels={$50$,$100$,$200$,$400$,$800$,$1600$,$2500$},
xticklabel style = {font=\small},
xminorticks=false,
ymin=1e-13,
ymax=1,
yminorticks=true,
ylabel style={font=\small\color{white!15!black}},
ylabel={Relative $\ell^{2}$ error},
ytick={1,1e-1,1e-2,1e-3,1e-4,1e-5,1e-6,1e-7,1e-8,1e-9,1e-10,1e-11,1e-12,1e-13,1e-14},
yticklabels={$10^{0}$,$$,$10^{-2}$,$$,$10^{-4}$,$$,$10^{-6}$,$$,$10^{-8}$,$$,$10^{-10}$,$$,$10^{-12}$,$$,$10^{-14}$},
yticklabel style = {font=\small},
axis background/.style={fill=white},
xmajorgrids,
xminorgrids,
ymajorgrids,
legend style={font=\normalsize,at={(0.9,0.9)}, anchor=north east, legend cell align=left, align=left, draw=white!15!black},
clip mode=individual,
]
\addplot [color=mycolor1, line width=1.0pt, mark=o, mark options={solid, mycolor1}, mark size=2.0pt]
  table[row sep=crcr]{%
64.0 0.0005459605055104123\\
128.0 3.0812620250837056e-5\\
256.0 2.8955012268665306e-8\\
512.0 1.0275713861290834e-11\\
1024.0 1.1256759492891915e-12\\
2048.0 1.10264143824528e-12\\
};
\addlegendentry{RBF}

\addplot [color=mycolor2, line width=1.0pt, mark=+, mark options={solid, mycolor2}, mark size=2.0pt]
  table[row sep=crcr]{%
64.0 0.002440989906706032\\
128.0 6.057973328614493e-6\\
256.0 3.2024472294834597e-9\\
512.0 1.4273951603988204e-11\\
1024.0 1.1144757736711282e-12\\
2048.0 1.102709004471271e-12\\    
};
\addlegendentry{Rational}


\addplot [color=mycolor4, line width=1.0pt, mark=triangle, mark options={solid, mycolor4}, mark size=2.0pt]
  table[row sep=crcr]{%
64.0 0.00030550548408220203\\
128.0 5.95969552661293e-6\\
256.0 2.807870484224073e-9\\
512.0 1.4282503365102066e-11\\
1024.0 1.114446833031851e-12\\
2048.0 1.1027322990012208e-12\\
};
\addlegendentry{Analytic}

\addplot [color=black, line width=1.0pt, mark options={solid, black}]
  table[row sep=crcr]{%
64.0 0.892401328481747\\
128.0 0.0034859426893818243\\
256.0 1.3616963630397751e-5\\
512.0 5.3191264181241216e-8\\
1024.0 2.077783757079735e-10\\
2048.0 8.116342801092715e-13\\
};
\addlegendentry{$8$th}

\end{loglogaxis}
\end{tikzpicture}%


}}%
    \caption{(a) Extended function $f = \Delta u$ on a logarithmic scale, 
with $u$ from \cref{eq:harmext}, using seven levels of refinement. 
(b) and (c) are convergence plots for the relative $\ell^{\infty}$ error and the relative $\ell^{2}$ error, respectively, as the quad-tree is uniformly refined, with 
markers at refinement levels three to eight.}%
    \label{fig:rat}%
\end{figure}
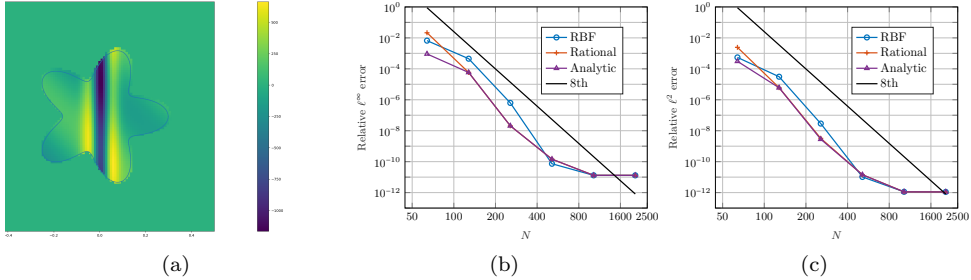

An alternate to our function extension scheme is to carry out
{\em one-dimensional} extension along lines in the plane. 
Consider a star-shaped domain centered at the origin, as shown in 
\cref{fig:rat}, on which we seek to solve the interior problem with 
solution \cref{eq:harmext}. 
For each point $\mathbf{x}$ outside $\Omega$ on the $\bar{K}\times \bar{K}$ 
grids 
for each $\bar{S}$, we extend along the line passing through the origin and 
$\mathbf{x}$.
Assuming $\mathbf{x} \in \bar{S}$ for some cut cell $S$ of side length $L$,
we assume we are given the data at eight uniformly-spaced interior nodes
over a distance $L$ from the boundary. We then form the one-dimensional 
barycentric rational interpolant with Floater-Hormann weights \cite{baryrat}, 
using the Julia implementation from \cite{baryratjl}. 
We then evaluate the interpolant at $\mathbf{x}$. 
The results are shown in \cref{fig:rat}. Note that, using this extension 
method, provides errors of about the same magnitude as the analytic extension.
Note also that we are {\em not} extending in the normal direction, but in the radial direction which intersects the boundary at some unspecified angle. 
The cost of this version of function extension is much less than that of
a VFMM call. When considering geometries such as \cref{fig:interiorsaw},
a more careful implementation will be required to avoid caustics.
We will return to this topic in the next section.

\section{Conclusions}
\label{sec:concl}

We have presented a potential theory-based solver for the Poisson equation 
in complicated two-dimensional geometries. To avoid computing a volume
potential over the actual domain, which involves 
complicated quadratures over cut leaf nodes in a quad-tree discretization,
we have developed a fast, high-order scheme to extend the source density 
smoothly to a slightly larger region where a volume integral
FMM (VFMM) can be applied \cite{vfmm}. The VFMM computes a volume
integral in linear time on an adaptive quad-tree, assuming that the 
source distribution is available on a tensor product grid for 
every leaf node in the tree. 
Unlike many earlier function extension schemes, 
we do not require the extended function to decay smoothly to zero. 
It is sufficient for it to be extended a sufficient distance from the
domain boundary (on the order of a single cut square width).
Having computed the volume potential, an auxiliary 
integral equation is solved to impose the desired boundary (and radiation)
condition.

The order of convergence of our scheme is dictated by the underlying 
discretization, not the extension method, since we can adjust the 
order of accuracy of our Gaussian interpolant to match that of the underlying
scheme. To make our extension efficient, we designed a
single, universal interpolation matrix which can be precomputed and
used for every cut leaf node which is intersected by the domain boundary.
From the universal interpolation matrix, rows are extracted corresponding
to data that lies in the domain interior. 
This leads to a small least squares problem that is solved
by $QR$ factorization on each cut square.
Unlike the earlier high-order extension method of
\cite{fryklund2018partition}, the present scheme visits each
extension square once, without blending, making it much simpler to 
implement in parallel. Furthermore, the extension scheme does not rely
on the smoothness of the boundary - just on being resolved by the 
user-provided data. The robustness, order of convergence, and accuracy 
of the scheme have been demonstrated with several numerical examples.

{For our interior problem,
with 11 digits of accuracy and eighth order convergence,
the VFMM itself runs at about 1M points/sec/core and 
the RBF-based function extension runs at about 500,000 
points/sec/core. 
The integral equation cost should be negligible
(it is linear scaling in the number of boundary points, but
sublinear in the total number of unknowns). It dominates here,
since we rely on a non-optimized iterative FMM-based scheme,
but the full solver still requires only about two seconds for a problem with
more than 200,000 unknowns.
We expect that with some modest modifications, the full solver should
achieve a throughput of close to 500,000 points/sec/core.}

A natural extension of the method presented here is to the three-dimensional
case. The main ingredients are available, such as 
high performance, parallelized VFMM libraries \cite{malhotra_biros_2015} and layer potential FMMs for boundary integral equations \cite{GREENGARD2021100092}. 
However, it remains to be determined how well the RBF-QR algorithm performs 
in three dimensions \cite{Fornberg2011StableCW}.
If the constants associated with the RBF-QR approach are too large, our
preliminary experiments, presented in \cref{ssec:extradial}, 
suggest that one-dimensonal extension may be equally effective and faster.
We have begun exploring the extension method of \cite{epsteinjiang2022},
which appears to be just as efficient as rational approximation, both in terms of speed and accuracy. We suspect that, for robustness, this should
always be done in the normal direction and we are actively investigating 
this approach. Finally, we should note that our function extension scheme
is unrelated to the governing PDE - it can be used with any 
potential-theoretic approach to boundary value problems in complicated 
domains when fast solvers like the VFMM are available.

\section*{Acknowledgments}
We would like to thank Ludvig af Klinteberg for a Laplace solver implemented in Julia 0.6 for smooth domains and Lukas Bystricky for an RCIP-based  Laplace solver 
implemented in \textsc{Matlab}{} for piecewise smooth domains.
We would also like to thank Samuel Potter, Charlie Epstein, Shidong Jiang, and Manas Rachh for several
helpful conversations.
\bibliographystyle{siamplain}
\bibliography{references}
\end{document}